\documentclass{amsart}
\usepackage{latexsym}

\usepackage{amssymb,amsmath,amscd}
\usepackage{bm}
\usepackage{amsfonts}
\usepackage{epsf}
\usepackage{graphicx}
\newtheorem{theorem}{Theorem}[section]
\newtheorem{lemma}[theorem]{Lemma}

\newtheorem{remark}[theorem]{Remark}

\theoremstyle{definition}

\numberwithin{equation}{section}

  \numberwithin{equation}{section}

\def\boxit#1{\vbox{\hrule height1pt \hbox{\vrule width1pt\kern1pt
     #1\kern1pt\vrule width1pt}\hrule height1pt }}

\newcommand{\bc}{\begin{center}}
\newcommand{\ec}{\end{center}}
\newcommand{\be}{\begin{eqnarray}}
\newcommand{\ee}{\end{eqnarray}}
\newcommand{\nn}{\nonumber}
\newcommand{\ben}{\begin{eqnarray*}}
\newcommand{\een}{\end{eqnarray*}}

\newcommand{\Om}{\Omega}

\newcommand{\pa}{\partial}
\newcommand{\si}{\sigma}

\def\x{\times}
\def\hK{\hat{K}}

\def\cE{\mathcal{E}}

\def\cT{\mathcal{T}}
\def\R{\mathbb{R}}
\def\S{\mathbb{S}}

\def\div{\operatorname{div}}

\DeclareMathOperator{\sspan}{span}
\DeclareMathOperator{\sym}{sym}

\begin{document}


\title[Conforming mixed
finite elements]{A new family of efficient conforming
mixed finite elements on both rectangular and cuboid meshes for linear  elasticity in the symmetric
formulation}

\author{Jun Hu}
\address{LMAM and School of Mathematical Sciences,
 Peking University, Beijing 100871, P. R. China\\ hujun@math.pku.edu.cn}

\date{\today}
\thanks{The  author was supported by  the NSFC Projects 11271035,  91430213 and 11421101.}

\begin{abstract}
 A new family of mixed finite elements is proposed for solving
the classical Hellinger--Reissner mixed problem of the elasticity equations.
  For two dimensions,  the normal stress  of  the matrix-valued stress field  is approximated by an enriched
Brezzi--Douglas--Fortin--Marini element of order $k$, and  the shear
stress by the serendipity element of order $k$,  the displacement field by an
enriched discontinuous vector-valued $P_{k-1}$ element. The  degrees of
freedom on each element of the lowest order element, which is of
first order,   is $10$ plus $4$.  For three dimensions, the normal stress is approximated
by  an  enriched Raviart--Thomas element of order $k$, and  each component of the shear stress  by a product space of
 the serendipity element space of two  variables and  the space of polynomials of degree $\leq k-1$
 with respect to  the rest variable, the displacement field by an enriched discontinuous vector-valued $Q_{k-1}$ element.
The  degrees of freedom on each element of the lowest order element, which is of
first order,  is $21$ plus $6$. A family of reduced elements is also proposed by dropping some interior bubble functions of the stress
and employing the discontinuous vector-valued $P_{k-1}$ (resp. $Q_{k-1}$) element for the displacement field on each element.
 As a result the lowest order elements have $8$ plus $2$ and $18$ plus $3$ degrees of freedom on each element for
 two and three dimensions,  respectively.

The  well-posedness condition and the optimal a priori error
estimate are proved for  this family of finite elements. Numerical
tests are presented to confirm the theoretical results.
\end{abstract}
\maketitle

\section{Introduction}
The first order system of equations,  for the symmetric stress field
   $\sigma\in\Sigma:=H({\rm div},\Omega,\mathbb {S})$
 and the displacement field $u\in V:=L^2(\Omega,\mathbb{R}^n)$, reads:
 Given $f\in L^2(\Om, \R^n)$
find $(\sigma,u)\in \Sigma\times V$  such that
\begin{equation}\label{continuous}
\begin{split}
&(A\sigma,\tau)_{L^2(\Om)}+(\div\tau,u)_{L^2(\Om)}=0,\\
&(\div\sigma,v)_{L^2(\Om)}=(f,v)_{L^2(\Om)}, \\
\end{split}
\end{equation}
for any $(\tau,v)\in \Sigma\times V $.
Here and throughout this paper, the compliance tensor
$A(x):\S\rightarrow \S$  is bounded and symmetric positive
definite uniformly for $x\in \Om$ with $\S:=\R^{n\x n}_{\sym}$
the set of symmetric tensors. The space $H(\div,\Om, \S)$ is
defined by
$$
H(\div,\Om,\S):=\{\tau\in L^2(\Om,\S) \text{ and },
 \div\tau\in L^2(\Om, \R^n)\}
$$
equipped with the norm
$$
\|\tau\|_{H(\div,\Om)}^2:=\|\tau\|_{L^2(\Om)}^2+
\|\div\tau\|_{L^2(\Om)}^2.
$$
The stress-displacement formulation within the
 Hellinger-Reissner  principle for the linear elasticity is one celebrated example
    of \eqref{continuous}.

Compared with the mixed formulation of the Poisson equation, see for
instance, \cite{Brezzi-Fortin}, there is an additional symmetric
requirement on the stress tensor.  Such a constraint makes the
stable discretization of the piecewise polynomials extremely
difficult.  Then one idea that may be come up with is to  enforce
the symmetry condition weakly, which in fact leads to Lagrange
multiplier methods  \cite{Amara-Thomas, Arnold-Brezzi-Douglas,
Boffi-Brezzi-Fortin, Morley, Stenberg-1, Stenberg-2, Stenberg-3}. As
an alternative method,  composite elements were proposed by Johnson and Mercier \cite{Johnson-Mercier},
and Arnold, Douglas  Jr., Gupta,  \cite{Arnold-Douglas-Gupta}.  That idea might be motivated by  the Hsieh-Clough-Tocher
element for the biharmonic problem \cite{CiaBook}.  Indeed, there is
an observation in \cite{Johnson-Mercier} that the discrete
divergence free space therein  is  the range of the Airy stress
function of the Hsieh-Clough- Tocher plate element space, see a similar observation in \cite{Arnold-Douglas-Gupta}.  Given  a
scalar field $q$,  the Airy stress function  reads
$$
Jq:=
\begin{pmatrix}
\frac{\pa^2 q}{\pa y^2}& -\frac{\pa^2 q}{\pa x\pa y}\\
-\frac{\pa^2 q}{\pa x\pa y} & \frac{\pa^2 q}{\pa x^2}
\end{pmatrix}.
$$
Unfortunately,  this observation  was not further explored until
more than twenty years later its importance was realized
  by Arnold and Winther \cite{Arnold-Winther-conforming}.
  In that landmark   paper, it was  found that to
  design a stable discrete scheme is to look for a discrete
   differential complex with  the commuting diagram which reads, for two dimensions,
  \[
\begin{CD}
0 @>>> P_1(\Om) @>\subset>> C^{\infty}(\Om) @>J>>
C^{\infty}(\Om,\S) @>\div>> C^{\infty}(\Om,\R^2)  @>>> 0\\
&&@VVidV  @VVI_h V @VV\Pi_h V @VVP_hV\\
0 @>>> P_1(\Om) @>\subset>> Q_h @>J_h>> \Sigma_h   @>\div_h>>
V_h @>>> 0
\end{CD}
\]
where $Q_h$ is some conforming or nonconforming finite element space
for the  biharmonic equation; $J_h$ and $\div_h$ are the discrete
counterparts of the Airy operator $J$ and the divergence operator
$\div$, respectively, with respect to some regular triangulation
$\mathcal{T}_h$ of $\Omega$;
 $\Sigma_h$ and $V_h$ are some finite element approximations of $\Sigma$ and $V$, respectively;
  $I_h$ and $\Pi_h$ are canonical interpolation
operators  for the spaces $Q_h$ and $\Sigma_h$, respectively; $P_h$
is the $L^2$ projection operator from $V$ onto $V_h$.  In
particular, this commuting diagram implies  the Fortin Lemma
\cite{Brezzi-Fortin}. See,   Arnold, Awanou, and Winther
\cite{Arnold-Awanou-Winther} for the corresponding theory in three
dimensions.  Based on those fundamental theories,  conforming mixed
finite elements of piecewise polynomials  on both simplicial and
product meshes can  then  be  developed for both 2D and
3D \cite{Adams-Cockburn,Arnold-Awanou,Arnold-Awanou-Winther,
Arnold-Winther-conforming}; see
\cite{Carstensen,Carstensen-Gunther-Reininghaus-Thiele2008} for the
implementation of the lowest order method of
\cite{Arnold-Winther-conforming}. To avoid  complexity of
conforming mixed elements, several remedies are proposed, see,
\cite{Arnold-Falk-Winther, Cockburn-Gopalakrishnan-Guzman,
Gopalakrishnan-Guzman, Guzman} for new weak-symmetry finite
elements, \cite{Arnold-Winther-n,Gopalakrishnan-Guzman-n,Hu-Shi,
Man-Hu-Shi,Yi} for non-conforming finite elements. See also
\cite{Awanou,Chen-Wang} for the enrichment of  nonconforming
elements of \cite{Hu-Shi,Man-Hu-Shi} to conforming elements.  In a
recent paper \cite{HuManZhang2012}, a family of first order
 nonconforming mixed finite elements on product meshes  is proposed  for the
first order system of equations in any dimension, which
was extended to a family of conforming mixed elements in \cite{HuManZhang2013}.

This paper presents a family of conforming  mixed
elements for both two and three dimensions ($n=2, 3$), which can be regarded as a generalization to any order of the first order methods
  from \cite{HuManZhang2013}. It is motivated by an observation that the conformity of
the discrete  methods on  product  meshes  can be guaranteed by
the $H(\div)$-conformity of the  normal stress and the
$H^1$-conformity of two corresponding variables for each  component
of the shear stress; see also \cite{BecacheJolyTsogka2002} and \cite{HuManZhang2013,HuManZhang2012} for a similar observation in two dimensions.
For two dimensions, in these elements,  an enriched
 Brezzi-Douglas-Fortin-Marini (hereafter BDFM) element of
order $k$ is proposed to approximate the normal stress,
 the  serendipity element of order $k$ \cite{Arnold-Awanou2011,BrennerScott,CiaBook} is used to approximate the shear stress.
This discrete  space for the stress and an enriched discontinuous $P_{k-1}$ element  for the  displacement space
are able to form a stable discretization of the two dimensional problem under consideration.
 In the first order method which is the two dimensional element of \cite{HuManZhang2013},  of this family,  the total degrees
of freedom is,   $|E|+6|K|$ +$|P|$,  with $|E|$ the number of
edges, and $|K|$ the number of elements, and $|P|$ the number of
vertices of the partition $\cT_h$. Note that the total degrees of freedom of
 the first order conforming mixed  element method on rectangular meshes in  \cite{Chen-Wang} is, $3|E|+6|K|$ +$|P|$.
For three dimensions, an enriched Raviart--Thomas element of order $k$ is constructed to approximate the normal stress,
 and each component of the shear stress is approximated by a product space of the serendipity element  of order $k$  with respect to
 two associated variables and the $P_{k-1}$ element with respect to the rest variable. An enriched  $Q_{k-1}$ element space
  is taken as the space for the displacement.   In the first order method which is the three dimensional element of \cite{HuManZhang2013},  of this family,  the total degrees of freedom is,   $|E|+|F|+9|K|$ ,  with $|E|$ the number of
edges, $|F|$ the number of faces, and $|K|$ the number of elements,  of the partition $\cT_h$.
Note that the total degrees of freedom of  the first order conforming mixed element method on  cuboid  meshes in  \cite{Awanou} is,
$2|E|+8|F|+18|K|$.  A family of reduced elements is also proposed by dropping  interior bubble functions
on each element. As a result the lowest order elements have $8$ plus $2$ and $18$ plus $3$ degrees of freedom on each element for
 two and three dimensions,  respectively, which were  announced independently in  \cite{Chen2013} after the first version of this paper was submitted.

These spaces of this paper  are perfectly and tightly matched on
each element. However, the analysis of the discrete inf-sup
conditions for these elements has to overcome the difficulty of not
using directly the Fortin Lemma, the key ingredient for the stability
analysis of the mixed finite element method for the elasticity
problem, see, for instance,
\cite{Adams-Cockburn,Arnold-Awanou,Arnold-Awanou-Winther,
Arnold-Winther-conforming}.  For pure displacement boundary  problem,  the remedy  is  an explicitly
constructive proof of the discrete inf-sup condition, which can be
regarded as a generalization to the more general case of the idea
due to \cite{HuManZhang2012}; see also \cite{BecacheJolyTsogka2002} and \cite{HuManZhang2013}.
For the more general case,  in particular the pure traction boundary  problem,  we prove that the divergence space of the
$H(\div)$ bubble function space  is identical to the orthogonal complement space of the rigid motion space
with respect to the discrete displacement space on each  macro-element.  As we shall see in Section 4,  the proof for
 such a result is  very difficult and   complicated. One important technique  is to use two classes of orthogonal polynomials, namely, the Jacobi polynomials and the Legendre polynomials.  As a second step, we construct a quasi--interpolation operator
 to control macroelementwise  rigid motion for $k>1$. Then the discrete inf--sup condition follows.  For the first order methods with $k=1$,
we succeed in proposing a new macroelement technique to  finally establish the
discrete inf--sup condition, which can be regarded  as an extension to the more  general case of that from \cite{Stenberg-1,Stenberg90}.

This paper is organized as follows.  In  the following two sections, we
 present the new mixed elements for two dimensions and analyze their  properties
including the well-posedness. In section 4, we  consider the pure traction boundary  problem and prove
the well-posedness of the discrete problem. In section 5 we define the new  mixed elements for three dimensions.
In section 6, we present a family of reduced elements by dropping some interior bubble functions
on each element. In section 7 we briefly summarize the error estimates of the discrete solutions and present two numerical examples,
 one for the pure  displacement boundary  problem, and the other for the pure traction boundary  problem.

\section{Mixed finite  element approximation  in two dimensions}
For approximating Problem \eqref{continuous} by the finite element
method, we introduce a rectangular triangulation $\cT_h$ of the
rectangular domain $\Om\subset \R^2$ such that
$\bigcup_{K\in \cT_h}K=\bar{\Omega}$, two distinct elements $K$ and
$K'$ in $\cT_h$ are either disjoint, or share the common edge $e$,
or a common vertex. Let $\cE$ denote the set of all edges in $\cT_h$
with $\cE_{K,V}$ the two vertical edges of $K\in\cT_h$ and
$\cE_{K,H}$ the two horizontal edges.  Let $\cE_V^I$ and $\cE_H^I$ denote the sets of all the interior vertical and horizontal edges of
 $\cT_h$, respectively, and  $\mathcal{V}^I$ be the  set of all the internal vertices of $\cT_h$.  Given  vertex $A\in \mathcal{V}^I$,
  let $\cE(A)$ be the set of edges that take $A$ as  one of their endpoints. Given any edge $e\in\cE$ we
assign one fixed unit normal $\nu$ with $(\nu_1,\, \nu_2)$ its components,
also let $t=(-\nu_2, \nu_1)$ denote the tangential vector.

For each $K\in \cT_h$,  we introduce the following affine invertible
transformation
$$
F_K:\hK\rightarrow K, x=\frac{h_{x,K}}{2}\xi+x_{0,K},
\quad y=\frac{h_{y,K}}{2}\eta+y_{0,K}
$$
with the center $(x_{0,K}, y_{0,K})$,  the horizontal and
vertical edge lengthes  $h_{x,K}$ and $h_{y,K}$,  respectively,
and the reference element $\hK=[-1,1]^2$.
Given any integer $k$, let $P_k(\omega)$ denote the space of
polynomials over $\omega$ of total degrees not greater than $k$, let $Q_k(\omega)$ denote the
 space of polynomials of degree not greater than $k$ in each variable.
 Let $P_k(X)$ be the space of polynomials of degree not greater than $ k$ with respect to the  variable $X$,
 and $P_k(X,Y)$ be the space of polynomials of degree not greater than $ k$ with  respect to the variables $X$ and $Y$.

 For the symmetric
fields $\sigma=\begin{pmatrix}
\sigma_{11} &\sigma_{12}\\
\sigma_{12} &\sigma_{22}
\end{pmatrix}\in\S$, we refer to $\sigma_n:=(\sigma_{11}, \sigma_{22})^T$ as the normal stress and $\sigma_{12}$ as the shear stress.

Before defining the space for the stress, we introduce  new mixed finite elements for the second order Poisson equation
 and the serendipity element of \cite{Arnold-Awanou2011,BrennerScott,CiaBook}.  Given
$K\in\mathcal{T}_h$ and an integer $k\geq 1$, the new mixed finite element space of order $k$ for the second order Poisson equation
reads :
\[
H_k(K):=(P_k(K))^2\backslash\ \sspan\{(0,x^k)^T,(y^k, 0)^T\}\oplus E_k(K),
\]
where
$$
E_k(K):=\sspan\{(x^{k+1}, 0)^T, (0, y^{k+1})^T, (x^2y^{k-1}, 0)^T, (0, y^2x^{k-1})^T\}.
$$

To define the degrees of freedom of the space $H_k(K)$, we introduce the well--known Jacobi polynomials:
\begin{equation}\label{Jacobi}
J_\ell(\xi):=((\ell+1)!)^2\sum\limits_{s=0}^\ell\frac{1}{s!(\ell+1-s)!(s+1)!(\ell-s)!}\bigg(\frac{\xi-1}{2}\bigg)^{\ell-s}\bigg(\frac{\xi+1}{2}\bigg)^s,
\end{equation}
 for any $\xi\in[-1, 1]$. The Jacobi polynomials satisfy the orthogonality condition:
\begin{equation}\label{orth}
\int_{-1}^1(1-\xi^2)J_l(\xi) J_m(\xi) d \xi=\frac{8}{2l+3}\frac{((l+2)!)^2}{(l+3)!l!}\delta_{lm}.
\end{equation}
We also need the Legendre polynomials
$$
L_\ell(\xi):=\frac{1}{2^\ell \ell!}\frac{d^\ell(\xi^2-1)^\ell}{d\xi^\ell} \text{ for any }\xi\in [-1, 1].
$$
The Legendre polynomials satisfy the orthogonality condition:
\begin{equation}\label{Lorth}
\int_{-1}^1L_l(\xi) L_m(\xi) d \xi=\frac{2}{2l+1}\delta_{lm}.
\end{equation}

\begin{lemma}\label{lemma2.1}
The vector-valued function $(\hat{q}_1, \hat{q}_2)^T=:\hat{q}\in H_k(\hat{K})$  can be uniquely
determined by the following conditions:
\begin{enumerate}
\item
$\int_{\hat{e}}\hat{q}\cdot \hat{\nu} \hat{p}d\hat{s}\text{ for any }\hat{p}\in P_{k-1}(\hat{e})
\text{ and  any }\hat{e}\subset \partial \hat{K}$,\\
\item
$\int_{\hat{K}} \hat{q}_1 J_{k-1}(\xi)d\xi d\eta$, \text{ and } $\int_{\hat{K}} \hat{q}_2 J_{k-1}(\eta)d\xi d\eta$,\\
\item
$\int_{\hat{K}} \hat{q}_1 L_{k-1}(\eta)d\xi d\eta$,\text{ and } $\int_{\hat{K}} \hat{q}_2 L_{k-1}(\xi)d\xi d\eta$,\\
\item
$\int_{\hat{K}} \hat{q}\cdot \hat{p} d\xi d\eta  \text{ for any }\hat{p}\in (P_{k-2}(\hat{K}))^2$.
\end{enumerate}
\end{lemma}
\begin{proof}  Since the  dimension of  the space $H_k(\hat{K})$ is equal to the number of these conditions,  it suffices to
 prove that $\hat{q}\equiv 0$ if  these conditions vanish. Since $\hat{q}\cdot \hat{\nu}\in P_{k-1}(\hat{e})$, the first  condition (1)
  implies that
  $$
  \hat{q}_1=(1-\xi^2)(\hat{g}_1+c_1J_{k-1}(\xi)+b_1L_{k-1}(\eta)), \text{ and } \hat{q_2}=(1-\eta^2)(\hat{g}_2+c_2J_{k-1}(\eta)+b_2L_{k-1}(\xi)),
  $$
  where $\hat{g}_1\,, \hat{g}_2\in P_{k-2}(\hat{K})$, and $c_1\,, c_2\,, b_1\,,b_2$ are four interpolation parameters, and $J_{k-1}$ and
  $L_{k-1}$ are the  Jacobi and  Legendre polynomials of degree $k-1$, respectively.  We first consider the case $k\geq 2$. It follows from \eqref{orth} that
   $$
   \int_{\hat{K}}(1-\xi^2)(\hat{g}_1+ b_1L_{k-1}(\eta))J_{k-1}(\xi)d\xi d\eta=0
   $$
    and
    $$
    \int_{\hat{K}}(1-\eta^2)(\hat{g}_2+b_2L_{k-1}(\xi))J_{k-1}(\eta) d\xi d\eta=0.
   $$
      Therefore, by the condition (2),
   $$
   c_1=c_2=0.
   $$
   The condition  \eqref{Lorth} implies
   $$
   \int_{\hat{K}}(1-\xi^2)\hat{g}_1L_{k-1}(\eta)d\xi d\eta=0
   $$
    and
    $$
    \int_{\hat{K}}(1-\eta^2)\hat{g}_2L_{k-1}(\xi)) d\xi d\eta=0.
   $$
   This and the condition (3) yield
   $$
   b_1=b_2=0.
   $$
   Hence the final result follows from the condition (4).  For the case $k=1$, the condition (2) is identical to the condition (3).
   A similar argument above completes the proof.
\end{proof}

\begin{remark}
The  space $H_k(K)$ is an enrichment of  the  BDFM
element space from \cite{Brezzi-Douglas-Fortin-Marini}.  Hence we call this new mixed element as the enriched
BDFM element.
\end{remark}

The global space of the enriched BDFM element reads
\[
H_k(\mathcal{T}_h):=\{q\in H(\div, \Omega, \R^2), q|_K \in
H_k(K) \text{ for any }K\in\cT_h \}.
\]
Note that, for any $q\in H_k(\mathcal{T}_h)$, the first component
of $q$ is continuous across the interior vertical edges of $\cT_h$
while the second component of $q$ is continuous across  the interior
horizontal edges of $\cT_h$.

 To get a stable pair of spaces,
 we propose to use the  serendipity element of order $k$ from \cite{Arnold-Awanou2011,BrennerScott,CiaBook} to approximate the
 shear  stress, which reads
  \[
  S_k(x,y):=P_k(x,y)+\sspan\{x^ky, xy^k\}.
  \]
Given any $\tau_{12}\in S_k(x,y)$, it can be uniquely determined
by the following conditions \cite{Arnold-Awanou2011}:
  \begin{enumerate}
  \item  the values of $\tau_{12}$ at four vertices of $K$,\\
\item the values of $\tau_{12}$ at $k-1$ distinct points in the interior of each edge of
$K$,\\
  \item  the moments $\int_K \tau_{12} pdxdy$ for any $p\in P_{k-4}(K)$.
  \end{enumerate}
  The global space of the serendipity element of order $k$ is
  defined as
 \[
 S_k(\cT_h):=\{\tau_{12}\in H^1(\Omega), \tau_{12}|_K\in S_k(x,y)\text{ for any }K\in\cT_h\}.
 \]
Note that the space $S_1(\cT_h)$ is the usual $H^1$-conforming bilinear element space.

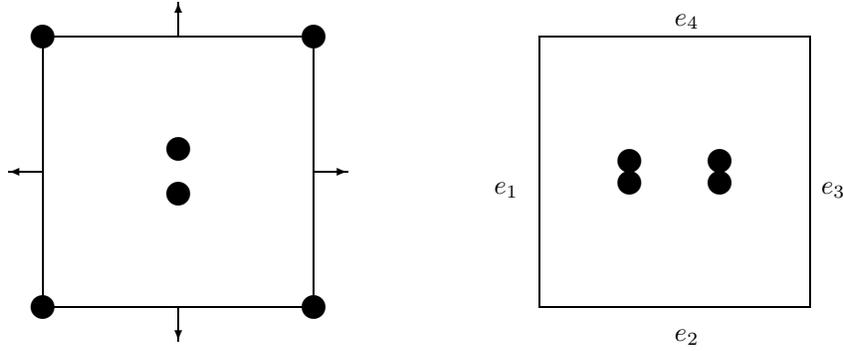
\begin{figure}[h]
\setlength{\unitlength}{0.30cm}
\begin{picture}(30,12)
\multiput (0,0)(0,12){2}{\line(1,0){12}} \multiput
(0,0)(12,0){2}{\line(0,1){12}}

\put(12,6){\vector(1,0){1.5}}

\put(0,6){\vector(-1,0){1.5}}

\put(6,12){\vector(0,1){1.5}}

\put(6,0){\vector(0,-1){1.5}}

\put(6,5){\circle*{1}} \put(6,7){\circle*{1}}

\multiput (0,0)(0,12){2} {\circle*{1}}
 \multiput(12,0)(0,12){2}{\circle*{1}}
\multiput (22,0)(0,12){2}{\line(1,0){12}} \multiput
(22,0)(12,0){2}{\line(0,1){12}}
\put(26,5.5){\circle*{1}}
 \put(26,6.5){\circle*{1}}
\put(30,5.5){\circle*{1}}
\put(30,6.5){\circle*{1}}
 \put(20,5){$e_1$} \put(34.5,5){$e_3$}
\put(28,-1.5){$e_2$} \put(28,12.5){$e_4$}
\end{picture}
\caption{Element diagram for the lowest order stress and
displacement}\label{Degrees}
\end{figure}

The discrete space of the element  is  combined from the enriched
BDFM element space and the serendipity element space:
\[
\Sigma_k(K):=\{\tau\in\S, \tau_n\in H_k(K), \tau_{12}\in S_k(x,y)\}.
\]
The degrees of freedom are inherited from the enriched  BDFM element and the
serendipity element:
\begin{enumerate}
\item  the moments of degree not greater than $k-1$ on the four edges of $K$
       for $\sigma_n\cdot \nu $,\\
\item the moments of  degree not greater than $k-2$ on  $K$ for $\sigma_n$,\\

\item the values
$\int_K (\sigma_n)_1 J_{k-1}(2(x-x_{0, K})/h_{x, K})dx dy$,\text{ and } $\int_K (\sigma_n)_2 J_{k-1}(2(y-y_{0, K})/h_{y, K})dx dy$
 where $(\sigma_n)_1$ is the first component of $\sigma_n$, and $(\sigma_n)_2$ is the second component of
  $\sigma_n$,\\

\item the values
$\int_K (\sigma_n)_1 L_{k-1}(2(y-y_{0, K})/h_{y, K})dx dy$,\text{ and } $\int_K (\sigma_n)_2 L_{k-1}(2(x-x_{0, K})/h_{x, K})dx dy$,\\

\item  the values of $\sigma_{12}$ at four vertices of $K$,\\
\item the values of $\sigma_{12}$ at $k-1$ distinct points in the interior of each edge of
$K$,\\
\item the moments of degree not greater than $k-4$ on  $K$ for
$\sigma_{12}$.
\end{enumerate}
The definitions of the enriched BDFM element and the serendipity element
imply that these conditions are unisolvent for the space
$\Sigma_k(K)$.  The degrees of freedom for the lowest order element
is illustrated in Figure \ref{Degrees}.

The global space of order $k$ is defined as
\begin{equation}\label{Sigmah}
\Sigma_k(\cT_h):=\{\tau\in \Sigma, \tau|_K\in \Sigma_k(K) \text{ for any }K\in\cT_h\}.
\end{equation}
On each element $K$, the space for the displacement is taken as
$$
V_k(K):=(P_{k-1}(K))^2\oplus\sspan\{(x^k, 0)^T, (0, y^k)^T, (xy^{k-1}, 0)^T, (0, x^{k-1}y)^T\}.
$$

Then the global space for the displacement  reads
\begin{equation}\label{Vh}
V_{k}(\cT_h):=\{v\in V, v|_K\in V_k(K) \text{ for any } K\in \cT_h\}.
\end{equation}
\begin{remark} The lowest order element (k=1) of this family has
 10 stress  and  4 displacement degrees of freedom  per element,
 which is the two dimensional element of \cite{HuManZhang2013}, see  degrees of freedom in Figure \ref{Degrees}.
\end{remark}

It follows from the definitions of  the spaces $\Sigma_k(\cT_h)$ and
$V_k(\cT_h)$ that $\div \Sigma_k(\cT_h)\subset V_k(\cT_h)$; in the
following section, we shall prove  the converse $V_k(\cT_h)\subset
\div \Sigma_k(\cT_h)$. This indicates the well-posedness of  this
family of elements.

The mixed element methods can be stated as: Find
$(\sigma_{k,h},u_{k,h})\in\Sigma_k(\cT_h)\x V_k(\cT_h)$ such that
\begin{equation}\label{discrete}
\begin{split}
&(A\sigma_{k,h},\tau)_{L^2(\Om)}
+(\div \tau,u_{k,h})_{L^2(\Om)}=0,\\
&(\div\sigma_{k,h},v)_{L^2(\Om)}=(f,v)_{L^2(\Om)},
\end{split}
\end{equation}
for any $(\tau,v) \in\Sigma_k(\cT_h)\x V_k(\cT_h)$.

\section{Well-posedness of discrete problem for pure displacement boundary  problem  in two dimensions}
In this section, we analyze the well-posedness of the discrete
problem \eqref{discrete}. From the mixed theory of
\cite{Brezzi-Fortin}, we need the following two assumptions
\begin{enumerate}
\item K-ellipticity. There exists a constant $C>0$
   independent of the meshsize  such that
\begin{equation}
 (A\tau,\tau)_{L^2(\Om)}\geq C\|\tau\|_{H(\div,\Om)}^2 \nn
\end{equation}
for any
\begin{equation}
\tau\in Z_k(\cT_h):=\{\tau\in\Sigma_k(\cT_h),  (\div\tau,
v)_{L^2(\Omega)}=0 \quad \text{ for all } v\in V_k(\cT_h)\}.\nn
\end{equation}
\item Discrete B-B condition. There exists a positive constant
  $C$ independent of the meshsize  with
\begin{equation}
\sup\limits_{0\not=\tau\in\Sigma_k(\cT_h)}\frac{(\div\tau,v)_{L^2(\Omega)}}
{\|\tau\|_{H(\div,\Om)}}\geq C\|v\|_{L^2(\Om)} \quad
 \text{ for any } v\in V_k(\cT_h).\nn
\end{equation}
\end{enumerate}
Herein and throughout, $C$ denotes a generic positive constant,
which may be different at the different occurrence but independent
of the meshsize $h$.  It follows from $\div\Sigma_k(K)\subset V_k(K)$ for
any $K\in\cT_h$ that $\div \tau=0$ for any $\tau\in Z_k(\cT_h)$.
This implies the  K-ellipticity condition.

To prove the discrete B-B condition, the usual idea in the
literature is to use the Fortin Lemma \cite{Brezzi-Fortin}.  More
precisely,  a bounded interpolation operator $\Pi_K:H^1(K,\S)
\rightarrow \Sigma_k(K)$ is constructed such that the following commuting
diagram property holds
\begin{equation}\label{commuting}
 \div\Pi_K\si=P_K\div\sigma \text{ for any } \sigma\in H^1(K,\S),
 \end{equation}
where $P_K$ is the projection operator from $L^2(K,\R^2)$ onto $V_k(K)$.
 So far, most of stable mixed  finite element methods for the
 linear elasticity problem within the Hellinger-Reissner principle
 are designed with such a property, see, for instance,
 \cite{Adams-Cockburn,Arnold-Awanou,Arnold-Awanou-Winther,
  Arnold-Winther-conforming}. However,  such a technique can not be used directly herein
   since there are not enough local degrees of freedom  for this family of elements under consideration.
  The idea  is to make a construction
  proof. More precisely,  given $v\in V_k(\cT_h)$, we  find explicitly
  $\tau\in \Sigma_k(\cT_h)$ such that
  \begin{equation}\label{eq3.1}
  \div \tau=v \ \ \ \text{ and } \|\tau\|_{H(\div, \Omega)}\leq
  C\|v\|_{L^2(\Omega)}.
  \end{equation}
  Such an idea is motivated by the stability analysis of the
  Raviart--Thomas element  for the Poisson equation in one dimension,
  which is first explored  to analyze the stability of a family of first order nonconforming
  mixed finite element methods on the product  mesh for the linear elasticity problem with the stress-displacement
  formulation in any dimension in a recent paper \cite{HuManZhang2012}.
  Therein,  the  discrete displacement is a piecewise  constant
  vector, which implies that the $\tau$ of  \eqref{eq3.1} can  be directly given so that $\div\tau=v$ for any $v$.
  In this paper we  use the form from \cite{HuManZhang2013}  to  construct $\tau$.

  For convenience,  suppose that the domain $\Omega$
     is a unit square $[0,1]^2$ which  is triangulated evenly into
    $N^2$ elements, $\{ K_{i j} \}$. This implies that
    $h_{x,K}=h_{y,K}=h:=1/N$ for any $K\in\cT_h$.
For any $v\in V_k(\cT_h)$, it can be decomposed as a sum,
$$
   v:=(v_1,v_2)^T=\sum\limits_{i=1}^{N}\sum\limits_{j=1}^Nv_{ij}\varphi_{ij}(x),
$$
where $\varphi_{ij}(x)$ is the characteristic function
 on the element $K_{ij}$ and $v_{ij}=(v_{ij}^1,v_{ij}^2)^T=v|_{K_{ij}}$.
Before the construction of $\tau\in\Sigma_k(\cT_h)$ with properties
of \eqref{eq3.1}, we need a
   decomposition of $v$.   We define the space $V_{y,ij}:=\sspan\{1,y,\cdots,
   y^{k-1}\}$, which introduces the following decomposition:
   \[
   P_{k-1}(K_{ij})\oplus\sspan\{x^{k},xy^{k-1}\}=V_{y,ij}\oplus (x-ih)\bigg(P_{k-2}(K_{ij})\oplus\sspan\{x^{k-1},y^{k-1}\}\bigg).
   \]
   This implies  that there exist unique  $v_{x,ij}^1\in P_{k-2}(K_{ij})\oplus\sspan\{x^{k-1},y^{k-1}\}$  and  $v_{y,ij}^1\in V_{y,ij}$
   such that
   \begin{equation}\label{decomp}
   v_{ij}^1=(x-ih)v_{x,ij}^1+v_{y,ij}^1.
   \end{equation}

\begin{theorem}\label{Theorem3.1}
It holds that
\begin{equation}
\sup\limits_{0\not=\tau\in\Sigma_k(\cT_h)}\frac{(\div\tau,v)_{L^2(\Omega)}}
{\|\tau\|_{H(\div,\Om)}}\geq \sqrt{\frac{2}{3}}\|v\|_{L^2(\Om)} \quad
 \text{ for any } v\in V_k(\cT_h).\nn
\end{equation}
\end{theorem}
\begin{proof} Given $v=(v_1, v_2)^T\in V_k(\cT_h)$, we define $T_{11}$ as the integration of $v_1$ along the rectangles along the $x$ direction:
   \begin{equation}\label{t11}
   \tau_{11}(x, y)=\int_0^x v_1(t,y)dt.
   \end{equation}
   On the element $K_{ij}:=[(i-1)h, ih]\times [(j-1)h, jh]$, $1\leq i, j\leq N$,
    by \eqref{decomp} and \eqref{t11}, it is  straightforward to see that
   $$
   \tau_{11}\in P_k(K) \backslash\ \sspan\{y^k\}\oplus\sspan\{x^{k+1},x^2y^{k-1}\},
   $$
  for $(x, y)\in K_{ij}$.  Similarly we can define $\tau_{22}$ as
  \begin{equation}\label{t22}
  \tau_{22}(x, y)=\int_0^yv_2(x,t)dt.
  \end{equation}
Then by \eqref{t11}, $\tau_{11}$ is continuous in the $x$ direction, by  \eqref{t22}, $\tau_{22}$ is continuous in the $y$
 direction. Hence we get an $H(\div)$ field
 $$
 \tau=\begin{pmatrix}\tau_{11}& 0\\ 0 &\tau_{22}\end{pmatrix}\in \Sigma_k(\cT_h).
 $$
 By the definition of $\tau$, it follows that
 \begin{equation}
 \div\tau=v.
 \end{equation}
 It remains to bound the $L^2$ norm of $\tau$.  We first consider the $L^2$ norm of the first component $\tau_{11}$:
\begin{equation*}
\begin{split}
\|\tau_{11}\|_{L^2(\Omega)}^2
&=\sum\limits_{i=1}^{N}\sum\limits_{j=1}^N\int_{K_{ij}}\bigg(\int_0^xv_1dt\bigg)^2dxdy\\
&\leq \sum\limits_{i=1}^{N}\sum\limits_{j=1}^N \int_{(j-1)h}^{jh}\int_{(i-1)h}^{ih}\bigg(x\int_0^xv_1^2dt\bigg)dx dy\\
&\leq \sum\limits_{i=1}^{N}\sum\limits_{j=1}^N \int_{(j-1)h}^{jh}\bigg(\int_0^{ih}v_1^2dt\bigg)\bigg( \int_{(i-1)h}^{ih} x dx \bigg) dy\\
&= \sum\limits_{i=1}^{N}\sum\limits_{j=1}^N\frac{h^2(2i-1)}{2} \int_{(j-1)h}^{jh}\int_0^{ih}v_1^2dtdy\\
&= h^2\sum\limits_{j=1}^{N}\sum\limits_{i=1}^N\sum\limits_{\ell=i}^N(\ell-\frac{1}{2})\|v_1\|_{L^2(K_{ij})}^2\\
&= h^2\frac{N(N-1)}{2}\sum\limits_{j=1}^{N}\sum\limits_{i=1}^N\|v_1\|_{L^2(K_{ij})}^2\leq \frac{1}{2} \|v_1\|^2_{L^2(\Omega)}.
\end{split}
\end{equation*}
A similar argument proves that
$$
\|\tau_{22}\|_{L^2(\Omega)}^2\leq \frac{1}{2} \|v_2\|^2_{L^2(\Omega)}.
$$
Hence
$$
\|\tau\|_{H(\div, \Omega)}^2=\|\div\tau\|_{L^2(\Omega)}^2+\|\tau\|_{L^2(\Omega)}^2\leq\frac{3}{2}\|v\|_{L^2(\Omega)}^2.
$$
This completes the proof.
\end{proof}

\section{The pure traction boundary  problem}
   This section considers the pure traction boundary  problem,   i.e.,  the stress space is subject to zero Neumann boundary
  condition while no boundary condition on the displacement.  In practice,  part of the elasticity body should be located, i.e,
    the displacement has a Dirichlet boundary condition on  some non-zero measure boundary.  But the pure traction boundary  problem
   is the most difficult one in mathematical analysis.  A similar  proof for Theorems \ref{Theorem4.1} and \ref{Theorem4.2}   can  prove them for partial displacement boundary  problems.

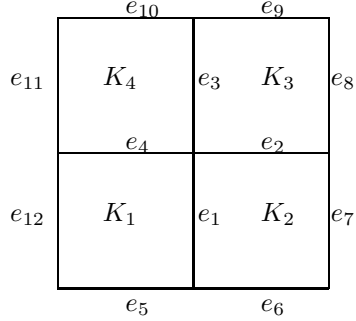
\begin{figure}[htb]
\begin{center}
\setlength{\unitlength}{0.30cm}
\begin{picture}(12,12)
\multiput (0,0)(0,6){3}{\line(1,0){12}}

 \multiput (0,0)(6,0){3}{\line(0,1){12}}

 \put (2,3){$K_1$}
 \put (9,3){$K_2$}
 \put (9,9){$K_3$}
\put (2,9){$K_4$}
\put (6.2,3){$e_1$}
\put (6.2,9){$e_3$}
\put (3,6.2){$e_4$}
\put (9,6.2){$e_2$}
\put (3, -1){$e_5$}
\put (9, -1){$e_6$}
\put (3, 12.2){$e_{10}$}
\put (9, 12.2){$e_9$}
\put (12.1,3){$e_7$}
\put (12.1,9){$e_8$}
\put (-2.1,3){$e_{12}$}
\put (-2.1,9){$e_{11}$}
\end{picture}
\caption{Macroelement}\label{Macroelement}
\end{center}
\end{figure}

Let $\text{RM }$ be the rigid motion space in  two dimensions, which reads
$$
\text{RM}:=\text{span}\bigg\{
\begin{pmatrix} 1\\ 0 \end{pmatrix}, \begin{pmatrix} 0\\ 1\end{pmatrix}, \begin{pmatrix}y\\ -x\end{pmatrix}\bigg\}.
$$
Consider a pure traction boundary  problem:

 \begin{equation}
 \label{e0}
\begin{split}
    \div \sigma  &= f  \quad \hbox{ in } \ \Omega:=(0,1)^2,\\
     \sigma    \nu &=0 \quad\hbox{ on } \ \partial\Omega, \\
     ( u, v) &=0  \quad \text{ for any }  v\in \text{RM},
     \end{split}
 \end{equation}
 where $\sigma:=A^{-1}\epsilon (u)$ for $u\in H^1(\Omega, \R^2)$.
By the same discretization of the uniform square grid $\cT_h$ with $h=1/N$
 as in the previous section,  the finite element equations  remain the same except
 the spaces are changed with boundary and rigid-motion free conditions:
\begin{equation}
\begin{split}
 (A\sigma_h, \tau)_{L^2(\Omega)} + (\div \tau, u_h)_{L^2(\Omega)} & = 0 \quad \text{ for all }\tau \in \Sigma_{k,0}(\cT_h), \\
     (\div \sigma_h, v)_{L^2(\Omega)} &=(f,v)_{L^2(\Omega)} \quad  \text{ for all }v \in V_{k,0}(\cT_h),
 \end{split}
 \end{equation}
where
   \begin{equation} \label{S-h-0}
   \begin{split}
    \Sigma_{k, 0}(\cT_h)
     &=\{ \tau=\begin{pmatrix}\tau_{11} & \tau_{12} \\ \tau_{12} & \tau_{22} \end{pmatrix}\in  \Sigma_k(\cT_h),
      \tau \nu=0 \ \quad  \text{ on  } \partial \Om\}, \\
	 V_{k,0}(\cT_h)
     &=\{ v=\begin{pmatrix}v_1 \\v_2 \end{pmatrix}\in V_k(\cT_h),
       \ (v, w)_{L^2(\Omega)}=0 \quad \text{ for all } w\in \text{RM}
       \}.
       \end{split}
     	\end{equation}
The earlier  analysis remains the same except the discrete B-B condition as the stress space $\Sigma_{k,0}(\cT_h)$ is  smaller than  $\Sigma_{k}(\cT_h)$.  To prove the discrete B-B condition for the pair   $(\Sigma_{k,0}(\cT_h),  V_{k,0}(\cT_h))$,  we introduce the concept of a macro-element, i.e., a union of four rectangles, see Figure \ref{Macroelement}.

Given a macroelement $M$, we define finite element spaces
$$
\Sigma_{k,0}(M):=\{\tau\in H(\div, \Om_M, \S), \tau|_K\in\Sigma_k(K) \text{ for any }K\subset  M, \tau\nu=0 \text{ on }\partial \Om_M \}
$$
and
$$
V_k(M):=\{v\in L^2(\Om_M, \R^2), v|_K\in V_k(K) \text{ for any }K\subset M\},
$$
where $\Om_M:=\bigcup\limits_{K\subset M}{\rm{int}(K)}$.  Define the orthogonal complement space of the rigid motion space $\text{RM}$ with  respect to
$V_k(M)$ by
$$
\text{RM}^\perp(M):=\{v\in V_k(M), (v, w)_{L^2(\Om_M)}=0 \text{ for any } w\in \text{RM}\}.
$$
\subsection{Discrete inf--sup conditions for higher order elements with $k\geq 2$}
In this subsection, we shall prove that, for $k\geq 2$,
$$
\div \Sigma_{k,0}(M)=\text{RM}^\perp(M),
$$
which helps to establish the discrete inf--sup conditions. To this end, we define the discrete kernel space of the divergence operator on the macroelement $M$ by
\begin{equation}\label{kernel}
N_M:=\{v\in V_k(M), (\div \tau, v)_{L^2(\Om_M)}=0\text{ for any }\tau\in  \Sigma_{k,0}(M)\}.
\end{equation}
We shall  show that  $N_M=\text{RM}$ for $k\geq 2$ to accomplish our goal.  The difficulty is how to explore  the local
    degrees of freedom for the shear stress.    One important technique is to invoke the Jacobi polynomials defined  in \eqref{Jacobi}, which needs the following polynomials:
    \begin{equation}
    \mathbb{J}_i(\xi)=\int_{-1}^\xi J_i(s)ds, i\geq 0, \xi\in [-1, 1].
    \end{equation}
For the four elements $K_i$, $i=1, \cdots, 4$, in the macroelement $M$ (see Figure \ref{Macroelement}),  we recall the  following affine mapping:
 $$
\xi_i=\frac{2x-2x_{0, K_i}}{h_{x,K_i}},
\quad \eta_i=\frac{2y-2y_{0, K_i}}{h_{y,K_i}}, (x,y)\in K_i,
$$
where  $(x_{0,K_i}, y_{0,K_i})$ is  the center of $K_i$,  $h_{x,K_i}$ and $h_{y, K_i}$ are  the horizontal and
vertical edge lengthes  of $K_i$,  respectively.  We also need the following spaces:
\begin{equation*}
\begin{split}
\Sigma_{12, e_1}:&=\bigg\{\tau_{12}\in H^1_0(\Om_M), \tau_{12}=q(1-\eta_1^2)\times\left\{\begin{array}{ll}(1+\xi_1) & \text{ on } K_1
                                                            \\ (1-\xi_2) & \text{ on } K_2 \end{array},  q\in P_{k-2}(\eta_1), \tau_{12}=0
                                                             \text{ on } K_3, K_4   \right.\bigg\},\\
\Sigma_{12, e_2}:&=\bigg\{\tau_{12}\in H^1_0(\Om_M), \tau_{12}=q(1-\xi_2^2)\times\left\{\begin{array}{ll}(1+\eta_2) & \text{ on } K_2
                                                            \\ (1-\eta_3) & \text{ on } K_3 \end{array},  q\in P_{k-2}(\xi_2), \tau_{12}=0
                                                             \text{ on } K_1, K_4   \right.\bigg\},\\
\Sigma_{12, e_3}:&=\bigg\{\tau_{12}\in H^1_0(\Om_M), \tau_{12}=q(1-\eta_3^2)\times\left\{\begin{array}{ll}(1-\xi_3) & \text{ on } K_3
                                                            \\ (1+\xi_4) & \text{ on } K_4 \end{array},  q\in P_{k-2}(\eta_3), \tau_{12}=0
                                                             \text{ on } K_1, K_2   \right.\bigg\}, \\
\Sigma_{12, e_4}:&=\bigg\{\tau_{12}\in H^1_0(\Om_M), \tau_{12}=q(1-\xi_4^2)\times\left\{\begin{array}{ll}(1-\eta_4)& \text{ on } K_4
                                                            \\ (1+\eta_1) & \text{ on } K_1 \end{array},  q\in P_{k-2}(\xi_4), \tau_{12}=0
                                                             \text{ on } K_2, K_3   \right.\bigg\}.
\end{split}
\end{equation*}
 The  restriction space on $M$ of the space $S_k(\cT_h)$ of  the serendipity element reads
\[
 S_{k, 0}(M):=\{\tau_{12}\in H^1_0(M), \tau_{12}|_{K_i}\in S_k(x,y), i=1, \cdots, 4\}.
 \]
 Note that $\Sigma_{12, e_i}\subset S_{k, 0}(M)$, $i=1,2,3,4$.

\begin{lemma}\label{Lemma4.1} For $k\geq 4$, suppose that $(v_1, v_2)^T\in V_k(M)$ is of the form
\begin{equation}
    v_1|_{K_{i}}=a_{-1,i}+\sum\limits_{\ell=0}^{k-2}a_{\ell}\mathbb{J}_\ell(\eta_i)\text{ and }
    v_2|_{K_{i}}=b_{-1,i}+\sum\limits_{\ell=0}^{k-2}b_{\ell}\mathbb{J}_\ell(\xi_i),
    \end{equation}
 with $a_{-1, 1}=a_{-1, 2}$, $a_{-1, 3}=a_{-1, 4}$, $b_{-1, 1}=b_{-1, 4}$, $b_{-1, 2}=b_{-1, 3}$, and that
 $$
 \int_{\Om_M}\frac{\partial \tau_{12}}{\partial y} v_1+\frac{\partial \tau_{12}}{\partial x}v_2 dxdy=0\text{ for any } \tau_{12}\in S_{k, 0}(M),
 $$
 then
 $$
 a_{-1,1}=a_{-1,2}=a_{-1,3}=a_{-1,4}, b_{-1,1}=b_{-1,2}=b_{-1,3}=b_{-1,4},
\frac{2a_0}{h_{y, k_i}}=-\frac{2b_0}{h_{x,k_i}}, a_{\ell}=b_{\ell}=0, \ell=1, \cdots, k-2.
$$
\end{lemma}
 \begin{proof} An integration by parts yields
\begin{equation} \label{eq4.7}
\begin{split}
0=\int_{\Om_M}\frac{\partial \tau_{12}}{\partial y} v_1+\frac{\partial \tau_{12}}{\partial x}v_2 dxdy
 &=-\sum\limits_{i=1}^4\int_{K_i}\tau_{12}\big(\frac{\partial v_1}{\partial y} +\frac{\partial v_2}{\partial x}\big) dxdy\\
 &\quad +\int_{e_1\cup e_3}\tau_{12}(a_{-1,1}-a_{-1,4})dx+\int_{e_2\cup e_4}\tau_{12}(b_{-1,1}-b_{-1,2})dy.
\end{split}
 \end{equation}
 We take $\tau_{12}$ in \eqref{eq4.7} such that
 \begin{equation*}
\tau_{12}|_{K_i}\in (1-\xi_i^2)(1-\eta_i^2)\sspan\{J_0(\eta_i), \cdots, J_{k-4}(\eta_i), J_1(\xi_i), \cdots, J_{k-4}(\xi_i)\}.
\end{equation*}
This  leads to
$$
\frac{2a_0}{h_{y, k_i}}=-\frac{2b_0}{h_{x,k_i}}, a_{\ell}=b_{\ell}=0, \ell=1, \cdots, k-4.
$$
To show these four parameters $a_{k-3}$, $a_{k-2}$, $b_{k-3}$ and $b_{k-2}$ to be zero, we  turn to the case where $k=4$. Since $J_0(\xi_i)=J_0(\eta_i)=1$, $J_1(\xi_i)=2\xi_i$ and $J_1(\eta_i)=2\eta_i$, we take $\tau_{12}\in \Sigma_{12, e_1}$ with $q=J_1(\eta_1)$ in \eqref{eq4.7}. Since $\int_{e_1}(1-\eta_1^2)(1+\xi_1)J_{1}(\eta_1) dy=0$, this yields
 $ a_1=0$.  Similarly, the choice of $\tau_{12}\in \Sigma_{12, e_4}$ with $q=J_1(\xi_1)$
shows $ b_1=0$. Then the choice of $\tau_{12}\in \Sigma_{12, e_1}$ with $q=J_0(\eta_1)$ in \eqref{eq4.7}, yields $a_{-1,1}=a_{-1,4}$;  while the choice
of $\tau_{12}\in \Sigma_{12, e_4}$ with $q=J_0(\xi_1)$ in \eqref{eq4.7}, leads to $b_{-1,1}=b_{-1,2}$. Hence we choose
$$
\tau_{12}\in \Sigma_{12, e_1}  \text
{ with } q=J_2(\eta_1),
$$
 $$
\tau_{12}\in \Sigma_{12, e_4}  \text
{ with } q=J_2(\xi_4),
$$
in \eqref{eq4.7}, respectively,  to show $a_2=b_2=0$.
Next  we consider the case where $k>4$ which allows to  take
$\tau_{12}\in \Sigma_{12, e_1}$ with $q=1$ in \eqref{eq4.7}.  This leads to $a_{-1,1}=a_{-1,4}$.  A similar argument with $\tau_{12}\in \Sigma_{12, e_4}$ and  $q=1$ gets $b_{-1,1}=b_{-1,2}$.
Therefore, the choices of
$$
\tau_{12}\in \Sigma_{12, e_1}  \text
{ with } q=J_{k-2}(\eta_1) \text{ and } q=J_{k-3}(\eta_1),
$$
 $$
\tau_{12}\in \Sigma_{12, e_4}  \text
{ with }  q=J_{k-2}(\xi_4) \text{ and } q=J_{k-3}(\xi_4),
$$
in \eqref{eq4.7}, respectively,  to prove
$$
a_{k-3}=b_{k-3}=a_{k-2}=b_{k-2}=0.
$$
This completes the proof.
 \end{proof}

\begin{lemma}\label{Lemma4.2} For $k=2, 3$, suppose that $(v_1, v_2)^T\in V_k(M)$ is of the form
\begin{equation}
    v_1|_{K_{i}}=a_{-1,i}+\sum\limits_{\ell=0}^{k-2}a_{\ell}\mathbb{J}_\ell(\eta_i)\text{ and }
    v_2|_{K_{i}}=b_{-1,i}+\sum\limits_{\ell=0}^{k-2}b_{\ell}\mathbb{J}_\ell(\xi_i),
    \end{equation}
 with $a_{-1, 1}=a_{-1, 2}$, $a_{-1, 3}=a_{-1, 4}$, $b_{-1, 1}=b_{-1, 4}$, $b_{-1, 2}=b_{-1, 3}$, and that
 $$
 \int_{\Om_M}\frac{\partial \tau_{12}}{\partial y} v_1+\frac{\partial \tau_{12}}{\partial x}v_2 dxdy=0\text{ for any } \tau_{12}\in S_{k, 0}(M),
 $$
 then
 $$
 a_{-1,1}=a_{-1,2}=a_{-1,3}=a_{-1,4}, b_{-1,1}=b_{-1,2}=b_{-1,3}=b_{-1,4},
\frac{2a_0}{h_{y, k_i}}=-\frac{2b_0}{h_{x,k_i}}, a_{\ell}=b_{\ell}=0, \ell=1, \cdots, k-2.
$$
\end{lemma}
\begin{proof} We only present the details for the case where $k=3$ since the proof for the case $k=2$ is similar and simple.
An integration by parts yields
\begin{equation} \label{eq4.7b}
\begin{split}
0=\int_{\Om_M}\frac{\partial \tau_{12}}{\partial y} v_1+\frac{\partial \tau_{12}}{\partial x}v_2 dxdy
 &=-\sum\limits_{i=1}^4\int_{K_i}\tau_{12}\big(\frac{\partial v_1}{\partial y} +\frac{\partial v_2}{\partial x}\big) dxdy\\
 &\quad +\int_{e_1\cup e_3}\tau_{12}(a_{-1,1}-a_{-1,4})dx+\int_{e_2\cup e_4}\tau_{12}(b_{-1,1}-b_{-1,2})dy.
\end{split}
 \end{equation}
 For such a case,  we have
 $$
 \frac{\partial v_1}{\partial y} +\frac{\partial v_2}{\partial x}|_{K_i}=\frac{2a_0}{h_{y, K_i}}+\frac{2b_0}{h_{x, K_i}}+\frac{2a_1J_1(\eta_i)}{h_{y,K_i}}
 +\frac{2b_1J_1(\xi_i)}{h_{x,K_i}}.
 $$
  Let ${\tau}_{12}\in \Sigma_{12, e_1}$ with $q=J_1(\eta_1)$ and ${\tau}_{12}\in \Sigma_{12, e_4}$ with $q=J_1(\xi_4)$ in \eqref{eq4.7b}, respectively. This yields  $a_1=0$ and $b_1=0$, respectively. The choices of    ${\tau}_{12}\in \Sigma_{12, e_1}$ with $q=J_0(\eta_1)$ and ${\tau}_{12}\in \Sigma_{12, e_4}$ with $q=J_0(\xi_4)$ yield,
  respectively,
  $$
  \frac{2a_0}{h_{y, K_i}}+\frac{2b_0}{h_{x, K_i}}+a_{-1,1}-a_{-1,4}=0 \text{ and }\frac{2a_0}{h_{y, K_i}}+\frac{2b_0}{h_{x, K_i}}+b_{-1,1}-b_{-1,2}=0.
  $$
 Now we let ${\tau}_{12}|_{K_1}=(1+\xi_1)(1+\eta_1)$ (with appropriate  definitions in   $K_2$, $K_3$, and $K_4$) in \eqref{eq4.7b} to obtain
 $$
 \frac{2a_0}{h_{y, K_i}}+\frac{2b_0}{h_{x, K_i}}+a_{-1,1}-a_{-1,4}+b_{-1,1}-b_{-1,2}=0
 $$
 Finally we solve these three equations to show  the desired result.
\end{proof}

\begin{lemma}\label{Lemma4.3} It holds, for $k\geq 2$,  that
\begin{equation}
\div \Sigma_{k,0}(M)=\text{RM}^\perp(M).
\end{equation}
\end{lemma}
\begin{proof}  Since it is   straightforward to see that $\div \Sigma_{k,0}(M)\subset \text{RM}^\perp(M)$, we only need to  prove that
$$
N_M=\text{span}\bigg\{
\begin{pmatrix} 1\\ 0 \end{pmatrix}, \begin{pmatrix} 0\\ 1\end{pmatrix}, \begin{pmatrix}y\\ -x\end{pmatrix}\bigg\}.
$$

Any  $v=(v_1,v_2)\in V_k(M)$ can be expressed as, for $\ell=1, \cdots, 4$,
$$
v_1|_{K_\ell}=\sum\limits_{i+j\leq k-1}a_{i,j}^{(\ell)}x^iy^j+a_{k,0}^{(\ell)}x^k+a_{1, k-1}^{(\ell)}xy^{k-1}
$$
and
$$
v_2|_{K_\ell}=\sum\limits_{i+j\leq k-1}b_{i,j}^{(\ell)}x^iy^j+b_{0, k}^{(\ell)}y^k+b_{k-1, 1}^{(\ell)}yx^{k-1}.
$$
We choose $\tau$ such that $\tau_{12}=\tau_{22}=0$ and
$$
\tau_{11}|_{K_{\ell}}\circ  F_{K_{\ell}}^{-1}\in (1-\xi^2)\times\hat{\Sigma}_{11, K_{\ell}}:=\sspan\{J_{k-1}(\xi), L_{k-1}(\eta)\}\oplus P_{k-2}(\hat{K}),  \ell=1, \cdots, 4.
$$
The condition for $N_M$ implies that
$$
0=\bigg(\frac{\partial \tau_{11}}{\partial x}, v_1\bigg)_{L^2(\Om_M)}=-\sum\limits_{\ell=1}^4\bigg(\tau_{11}|_{K_\ell}, \frac{\partial (v_1|_{K_\ell})}{\partial x}\bigg)_{L^2(\Om_M)}.
$$
Since $\frac{\partial (v_1|_{K_\ell})}{\partial x}\circ F_{K_{\ell}}^{-1}\in  \hat{\Sigma}_{11, K_{\ell}}$, this yields
$$
a_{k,0}^{(\ell)}=a_{1, k-1}^{(\ell)}=a_{i,j}^{(\ell)}=0 \text{ for all } i\geq 1 \text{ with } i+j\leq k-1.
$$
Hence $v_1$ is of the form
\begin{equation}\label{1}
v_1|_{K_\ell}=\sum\limits_{j\leq k-1}a_{0,j}^{(\ell)}y^j
\end{equation}
The continuity and degrees of $\tau_{11}$ across the edges $e_1$ and $e_3$  show (using the moments of degree not greater
 than $k-1$ on $e_1$ and $e_3$ for $\tau_{11}$)
\begin{equation}\label{2}
a_{0,j}^{(1)}=a_{0,j}^{(2)} \text{ and }a_{0,j}^{(3)}=a_{0,j}^{(4)},\quad j=0, \cdots, k-1.
\end{equation}
A similar argument for $v_{2}$ shows that $v_2$ is of the form
\begin{equation}\label{3}
v_2|_{K_\ell}=\sum\limits_{i\leq k-1}b_{i,0}^{(\ell)}x^i
\end{equation}
and
\begin{equation}\label{4}
b_{i,0}^{(1)}=b_{i,0}^{(4)} \text{ and }b_{i,0}^{(2)}=b_{i,0}^{(3)},\quad i=0, \cdots, k-1.
\end{equation}
To decide these parameters $a_{0,j}^{(\ell)}$ and $b_{i,0}^{(\ell)}$, we propose to use the degrees of freedom for the
 shear stress component $\tau_{12}$.  We choose $\tau$ such that $\tau_{11}=\tau_{22}=0$ and $\tau_{12}=\tau_{12}^{(\ell)}\in \Sigma_{12, e_\ell}$, $\ell=1, \cdots, 4$.   The condition for $N_M$, and \eqref{1}--\eqref{4}, produce

\begin{equation}\label{5}
\begin{split}
 0&= \int_{e_\ell}\tau_{12}^{(\ell)}[v_{(\ell+1 \text{ mod } 2)}]ds
-\int_{K_\ell}\tau_{12}^{(\ell)}\bigg(\frac{\partial v_1|_{K_\ell}}{\partial y}+\frac{\partial v_2|_{K_\ell}}{\partial x}\bigg)dxdy\\
&\quad -\int_{K_{(\ell+1 \text{ mod } 4)}}\tau_{12}^{(\ell)}\bigg(\frac{\partial v_1|_{K_{(\ell+1 \text{ mod } 4)}}}{\partial y}+\frac{\partial v_2|_{K_{(\ell+1 \text{ mod } 4)}}}{\partial x}\bigg)dxdy,
 \end{split}
\end{equation}
where $[\cdot]$ denotes the jump of piecewise functions across edge $e_\ell$.
Since $v_1$ (resp. $v_2$) is a piecewise polynomial with respect to variable $y$ (resp. $x$), the symmetries of $\tau_{12}^{(\ell)}$ with the edges $e_\ell$, $\ell=1, \cdots, 4$,  lead to
\begin{equation}\label{6}
\int_{K_1}\tau_{12}^{(1)}\frac{\partial v_1|_{K_1}}{\partial y}dxdy=\int_{K_2}\tau_{12}^{(1)}\frac{\partial v_1|_{K_2}}{\partial y}dxdy, \quad
\int_{K_3}\tau_{12}^{(3)}\frac{\partial v_1|_{K_3}}{\partial y}dxdy=\int_{K_4}\tau_{12}^{(3)}\frac{\partial v_1|_{K_4}}{\partial y}dxdy
\end{equation}
and
\begin{equation}\label{7}
\int_{K_2}\tau_{12}^{(2)}\frac{\partial v_2|_{K_2}}{\partial x}dxdy=\int_{K_3}\tau_{12}^{(2)}\frac{\partial v_2|_{K_3}}{\partial x}dxdy, \quad
\int_{K_4}\tau_{12}^{(4)}\frac{\partial v_2|_{K_4}}{\partial x}dxdy=\int_{K_1}\tau_{12}^{(4)}\frac{\partial v_2|_{K_1}}{\partial x}dxdy.
\end{equation}
Next let  all $\tau_{12}^{(\ell)}\in \Sigma_{12, e_\ell}$, $\ell=1, \cdots, 4$, be defined by, up to the variable $x$ or $y$, and some transformation(s),  the same polynomials $q$ of one variable of  degree $\leq k-2$. Since $[v_1]|_{e_2}=[v_1]|_{e_4}$ and $[v_2]|_{e_1}=[v_2]|_{e_3}$, the symmetries of $\tau_{12}^{(\ell)}$ imply  additionally  that
\begin{equation}\label{8}
\int_{e_1}\tau_{12}^{(1)}[v_2]ds=\int_{e_3}\tau_{12}^{(3)}[v_2]ds \text { and } \int_{e_2}\tau_{12}^{(2)}[v_1]ds=\int_{e_4}\tau_{12}^{(4)}[v_1]ds.
\end{equation}
A substitution of  equations \eqref{6} through \eqref{8} into \eqref{5}  shows  that
\begin{equation}\label{9}
\int_{K_1}\tau_{12}^{(1)}\frac{\partial v_1|_{K_1}}{\partial y}dxdy=\int_{K_2}\tau_{12}^{(1)}\frac{\partial v_1|_{K_2}}{\partial y}dxdy=
\int_{K_3}\tau_{12}^{(3)}\frac{\partial v_1|_{K_3}}{\partial y}dxdy=\int_{K_4}\tau_{12}^{(3)}\frac{\partial v_1|_{K_4}}{\partial y}dxdy
\end{equation}
and
\begin{equation}\label{10}
\int_{K_2}\tau_{12}^{(2)}\frac{\partial v_2|_{K_2}}{\partial x}dxdy=\int_{K_3}\tau_{12}^{(2)}\frac{\partial v_2|_{K_3}}{\partial x}dxdy=
\int_{K_4}\tau_{12}^{(4)}\frac{\partial v_2|_{K_4}}{\partial x}dxdy=\int_{K_1}\tau_{12}^{(4)}\frac{\partial v_2|_{K_1}}{\partial x}dxdy.
\end{equation}
Since both $\frac{\partial v_1|_{K_1\cup K_2}}{\partial y}$  and $\frac{\partial v_1|_{K_3\cup K_4}}{\partial y}$
are polynomials of degree $\leq k-2$ with respect to $y$,  the conditions of \eqref{9} for all $\tau_{12}^{(1)}\in \Sigma_{12, e_1}$ and
 $\tau_{12}^{(3)}\in \Sigma_{12, e_3}$, the conditions of \eqref{10} for all $\tau_{12}^{(2)}\in \Sigma_{12, e_2}$ and
 $\tau_{12}^{(4)}\in \Sigma_{12, e_4}$ show that $(v_1, v_2)^T\in V_k(M)$ is of the form
\begin{equation}
    v_1|_{K_{i}}=a_{-1,i}+\sum\limits_{\ell=0}^{k-2}a_{\ell}\mathbb{J}_\ell(\eta_i)\text{ and }
    v_2|_{K_{i}}=b_{-1,i}+\sum\limits_{\ell=0}^{k-2}b_{\ell}\mathbb{J}_\ell(\xi_i), i=1,\cdots,4,
    \end{equation}
 with $a_{-1, 1}=a_{-1, 2}$, $a_{-1, 3}=a_{-1, 4}$, $b_{-1, 1}=b_{-1, 4}$, $b_{-1, 2}=b_{-1, 3}$.
 Hence it follows from Lemmas \ref{Lemma4.1} and \ref{Lemma4.2} that
$$
 a_{-1,1}=a_{-1,2}=a_{-1,3}=a_{-1,4}, b_{-1,1}=b_{-1,2}=b_{-1,3}=b_{-1,4},
\frac{2a_0}{h_{y, k_i}}=-\frac{2b_0}{h_{x,k_i}}, a_{\ell}=b_{\ell}=0, \ell=1, \cdots, k-2.
$$
This completes the proof.
\end{proof}

\begin{lemma}\label{new} For any $v_h\in V_{k, 0}(\cT_h)$,  there  exists a $\tau_h\in \Sigma_{k, 0}(\cT_h)$ such that
\begin{equation}
\int_M (\div \tau_h-v_h)\cdot w dx=0 \text{ for any }  w\in \text{RM} \text{ and any macro--element } M
\end{equation}
and
\begin{equation}
\|\tau_h\|_{H(\div, \Omega)}\leq C\|v_h\|_{L^2(\Omega)}.
\end{equation}
\end{lemma}
\begin{proof} It is standard that there exists a $\tau:=\begin{pmatrix}\tau_{11} &\tau_{12}\\ \tau_{12} &\tau_{22}\end{pmatrix}\in H^1_0(\Omega, \S)$ such that
\begin{equation}
\div \tau=v_h \text{ and }\|\tau\|_{H^1(\Omega)}\leq C\|v_h\|_{L^2(\Omega)}.
\end{equation}
It follows from the degrees of freedom for the enriched BDFM element in Lemma \ref{lemma2.1} and for the serendipity element that  there exist $(\tau_{11,h}, \tau_{22,h})^T\in H_k(\cT_h)$ and $\tau_{12,h}\in S_k(\cT_h)$ such that,  for edges of macro-element $M$ (see Figure \ref{Macroelement} for notation),
\begin{equation*}
\int_{e_{7}\cup e_{8}} (\tau_{11}-\tau_{11,h})p dy=\int_{e_{11}\cup e_{12}} (\tau_{11}-\tau_{11,h})q dy=0 \text{ for any }p\in P_1(e_{7}\cup e_{8}), q\in P_1(e_{11}\cup e_{12}),
\end{equation*}
\begin{equation*}
\int_{e_{5}\cup e_{6}} (\tau_{22}-\tau_{22,h})p dx=\int_{e_{9}\cup e_{10}} (\tau_{22}-\tau_{22,h})q dx=0 \text{ for any } p\in P_1(e_{5}\cup e_{6}), q\in P_1(e_{9}\cup e_{10}),
\end{equation*}
\begin{equation*}
\int_{e_{11}\cup e_{12}} (\tau_{12}-\tau_{12,h}) dy=\int_{e_{7}\cup e_{8}} (\tau_{12}-\tau_{12,h}) dy=\int_{e_{5}\cup e_{6}} (\tau_{12}-\tau_{12,h}) dx=\int_{e_{9}\cup e_{10}} (\tau_{12}-\tau_{12,h})dx=0.
\end{equation*}
Let $\tau_h=\begin{pmatrix}\tau_{11,h} &\tau_{12,h}\\ \tau_{12, h} &\tau_{22, h}\end{pmatrix}$.  We additionally have
$$
\|\tau_h\|_{H(\div, \Omega)}\leq C\|\tau\|_{H^1(\Omega)}.
$$
This completes the proof.
\end{proof}

We are now ready to establish the following inf--sup condition.
\begin{theorem}\label{Theorem4.1} For $k\geq 2$, there exists a positive constant
  $C$ independent of the meshsize  with
\begin{equation}
\sup\limits_{0\not=\tau\in\Sigma_{k, 0}(\cT_h)}\frac{(\div\tau,v)_{L^2(\Omega)}}
{\|\tau\|_{H(\div,\Om)}}\geq C\|v\|_{L^2(\Om)} \quad
 \text{ for any } v\in V_{k, 0}(\cT_h).\nn
\end{equation}
\end{theorem}
\begin{proof} Given $v\in V_{k, 0}(\cT_h)$, it follows from Lemma \ref{new} that there exists a $\tau_1\in \Sigma_{k, 0}(\cT_h)$ such that
\begin{equation}
\int_M (\div \tau_1-v)\cdot w dx=0 \text{ for any }  w\in \text{RM} \text{ and any macro--element } M
\end{equation}
and
\begin{equation}
\|\tau_1\|_{H(\div, \Omega)}\leq C\|v\|_{L^2(\Omega)}.
\end{equation}
By Lemma \ref{Lemma4.3}, there exists a $\tau_2\in \Sigma_{k, 0}(\cT_h)$ such that
\begin{equation}
\div \tau_2=\div \tau_1-v \text{ and } \|\tau_2\|_{H(\div, \Omega)}\leq C\|\div \tau_1-v\|_{L^2(\Omega)}.
\end{equation}
Then we have $\div(\tau_1+\tau_2)=v$ and $\|\tau\|_{H(\div, \Omega)}\leq C\|v\|_{L^2(\Omega)}$.
\end{proof}

\subsection{Discrete inf--sup condition for the first order element with $k=1$} Since the analysis in the previous subsection
 can not be  applied to the current case,  it needs a separate analysis. The ingredient is a modified  macroelement technique.
 We also note that the  macroelement technique from \cite{Stenberg-1} can not be used  directly here since  the semi-norm $|\cdot|_{1,h,M}$ there
  is not equivalent to the semi--norm $|\cdot|_M$ there for the present case, see \cite[Theorem 4.1]{Stenberg-1}. To overcome this difficulty, for
  $v\in V_{1}(M)$,   we propose the following  mesh  dependent semi--norm, see Figure \ref{Macroelement} for notation,
   \begin{equation}\label{eq4.23}
   \begin{split}
   |v|_{1, h, M}^2&=\sum\limits_{i=1}^4\|\epsilon(v)\|_{0, K_i}^2+h_{e_1}^{-1}\|[v_1]\|_{L^2(e_1)}^2+h_{e_3}^{-1}\|[v_1]\|_{L^2(e_3)}^2\\
   &\quad +h_{e_2}^{-1}\|[v_2]\|_{L^2(e_2)}^2
   +h_{e_4}^{-1}\|[v_2]\|_{L^2(e_4)}^2\\
   &\quad+((v_1|_{K_1}-v_1|_{K_4})(\mathbb{M}(e_4))+(v_1|_{K_2}-v_1|_{K_3})(\mathbb{M}(e_2))\\
   &\quad +(v_2|_{K_1}-v_2|_{K_2})(\mathbb{M}(e_1))+(v_2|_{K_4}-v_2|_{K_3})(\mathbb{M}(e_3)))^2,
  \end{split}
   \end{equation}
   where $\mathbb{M}(e_i)$, $i=1, \cdots, 4$, denote the midpoints of edges $e_i$, and $[\cdot]$ denote the jump of piecewise functions over edge. 
   Define a global seminorm
   \begin{equation}
   |v|_{1,h}^2=\sum\limits_{M}|v|_{1, h, M}^2 \text{ for all  macro-elements consisting of four elements like that in Figure \ref{Macroelement}}.
   \end{equation}
   It is straightforward to see that $|\cdot|_{1, h}$ defines a norm over $V_{1, 0}(\cT_h)$.  For $\tau\in\Sigma_{1, 0}(\cT_h)$, we define the following
    mesh dependent norm:
    \begin{equation}\label{eq4.25}
    \|\tau\|_{0,h}^2=\|\tau\|_{L^2(\Omega)}^2+\sum\limits_{e\in\cE^I_V}h_e\|\tau_{11}\|_{L^2(e)}^2+\sum\limits_{e\in\cE^I_H}h_e\|\tau_{22}\|_{L^2(e)}^2
    +\sum\limits_{A\in\mathcal{V}^I}\sum\limits_{e\in\cE(A)}h_e^2\tau_{12}(A)^2.
    \end{equation}

\begin{lemma}\label{Lemma4.4} For any macroelement $M$ illustrated in Figure \ref{Macroelement}, it holds that
$$
N_M=\text{span}\bigg\{
\begin{pmatrix} 1\\ 0 \end{pmatrix}, \begin{pmatrix} 0\\ 1\end{pmatrix}, \epsilon_{y,-x}\bigg\},
$$
where
$$
\epsilon_{y,-x}:=\left\{ \begin{array}{ll}
\begin{pmatrix} -1\\ 1 \end{pmatrix} &\text{ on }K_1,\\
\begin{pmatrix} -1\\ -1 \end{pmatrix}& \text{ on }K_2,\\
\begin{pmatrix} 1\\ -1 \end{pmatrix} &\text{ on }K_3,\\
\begin{pmatrix}~~~ 1~~~\\~~~ 1~~~ \end{pmatrix} & \text{ on }K_4.
\end{array}\right.
$$
\end{lemma}
\begin{proof}
 Any  $v=(v_1,v_2)\in V_1(M)$ can be expressed as, for $\ell=1, \cdots, 4$,
$$
v_1|_{K_\ell}=a_0^{(\ell)}+a_1^{(\ell)}x
$$
and
$$
v_2|_{K_\ell}=b_0^{(\ell)}+b_1^{(\ell)}y.
$$
We choose $\tau$ such that $\tau_{12}=\tau_{22}=0$ and
$$
\tau_{11}|_{K_{\ell}}\circ  F_{K_{\ell}}^{-1}=1-\xi^2,  \ell=1, \cdots, 4.
$$
The condition for $N_M$ implies that
$$
a_1^{(\ell)}=0, \ell=1, \cdots, 4.
$$
Similarly,
$$
b_1^{(\ell)}=0, \ell=1, \cdots, 4.
$$
Hence we can use the degrees on the edges $e_1$ and $e_3$ for $\tau_{11}$ (using the moments of degree zero on $e_1$ and $e_3$ for $\tau_{11}$) to shows that
$$
a_0^{(1)}=a_0^{(2)}\text{ and }a_0^{(3)}=a_0^{(4)}.
$$
A similar argument proves
$$
b_0^{(1)}=b_0^{(4)}\text{ and }b_0^{(2)}=b_0^{(3)}.
$$
At the end we use the degree of $\tau_{12}$ at the interior vertex of  $M$ to complete the proof.
\end{proof}
We need  another seminorm for the space $V_1(M)$:
\begin{equation}
|v|_M=\sup\limits_{0\not=\tau\in\Sigma_{1,0}(M) }\frac{(\div \tau, v)_{L^2(\Omega_M)}}{\|\tau\|_{0, h, M}}
\end{equation}
It follows from Lemma \ref{Lemma4.4} that the seminorm $|\cdot|_{1, h, M}$ is equivalent to the seminorm $|\cdot|_M$.  This allows for following
   a similar argument of \cite{Stenberg-1} and the references therein to prove the discrete inf--sup condition.

\begin{theorem}\label{Theorem4.2}
There exists a positive constant
  $C$ independent of the meshsize  with
\begin{equation}
\sup\limits_{0\not=\tau\in\Sigma_{k, 0}(\cT_h)}\frac{(\div\tau,v)_{L^2(\Omega)}}
{\|\tau\|_{0, h}}\geq C|v|_{1, h} \quad
 \text{ for any } v\in V_{k, 0}(\cT_h).\nn
\end{equation}
\end{theorem}

\section{Mixed finite element for  three dimensions}
We define a family of  conforming mixed finite element methods in
 three dimensions  in this section.  To this end,  let $\cT_h$  be   a cuboid  triangulation of the
cuboid domain $\Om\subset \R^3$ such that $\bigcup_{K\in \cT_h}K=\bar{\Omega}$.  On element $K\in\cT_h$,  for $k\geq 1$,  we define an enriched  Raviart--Thomas element space by
 $$
 H_{k}(K)=(P_{k, k-1, k-1}(K)\oplus E_{k,x})\times (P_{k-1, k, k-1}(K)\oplus E_{k,y})\times (P_{k-1, k-1, k}(K)\oplus E_{k,z}),
 $$
 where
\begin{equation*}
\begin{split}
P_{k, k-1, k-1}(K)=P_k(x)\times P_{k-1}(y)\times P_{k-1}(z),\\[0.5ex]
P_{k-1, k, k-1}(K)=P_{k-1}(x)\times P_k(y)\times P_{k-1}(z),\\[0.5ex]
P_{k-1, k-1, k}(K)=P_{k-1}(x)\times P_{k-1}(y)\times P_k(z)
\end{split}
\end{equation*}
and
\begin{equation*}
\begin{split}
E_{k,x}=x^{k+1}(P_{k-1}(y)+P_{k-1}(z)),\\[0.5ex]
E_{k,y}=y^{k+1}(P_{k-1}(z)+P_{k-1}(x)),\\[0.5ex]
E_{k,z}=z^{k+1}(P_{k-1}(x)+P_{k-1}(y)).
\end{split}
\end{equation*}
To  construct  the degrees of freedom of the space $H_k(\hat{K})$, we define
\begin{equation*}
\begin{split}
\Psi_{k-1}(\hat{K}):&=P_{k-2,k-1,k-1}(\hat{K})\times P_{k-1,k-2,k-1}(\hat{K})\times P_{k-1,k-1,k-2}(\hat{K}),\\[0.5ex]
\mathcal{J}_{k-1}(\xi):&=J_{k-1}(\xi)(P_{k-1}(\eta)+P_{k-1}(\zeta)),\\[0.5ex]
\mathcal{J}_{k-1}(\eta):&=J_{k-1}(\eta)(P_{k-1}(\xi)+P_{k-1}(\zeta)),\\[0.5ex]
\mathcal{J}_{k-1}(\zeta):&=J_{k-1}(\zeta)(P_{k-1}(\xi)+P_{k-1}(\eta)),\\[0.5ex]
\end{split}
\end{equation*}
for any $(\xi, \eta, \zeta)\in \hat{K}:=[-1,1]^3$.  We  recall that $J_{k-1}(\xi)$ is the Jacobi polynomial of degree $k-1$ with respect to
$\xi$, and $L_{k-1}(\xi)$ is the Legendre polynomial of degree $k-1$ with respect to $\xi$.
\begin{lemma}\label{lemma5.1}
The vector-valued function $(\hat{q}_1, \hat{q}_2, \hat{q}_3)^T=:\hat{q}\in H_k(\hat{K})$  can be uniquely
determined by the following conditions:
\begin{enumerate}
\item
$\int_{\hat{e}}\hat{q}\cdot \hat{\nu} \hat{p}d\hat{s}\text{ for any }\hat{p}\in Q_{k-1}(\hat{e})
\text{ and  any }\hat{e}\subset \partial \hat{K}$,\\
\item
$\int_{\hat{K}} \hat{q}_1 \hat{p}d\xi d\eta d\zeta \text{ for any }\hat{p}\in \mathcal{J}_{k-1}(\xi)$,\\

\item
$\int_{\hat{K}} \hat{q}_2 \hat{p}d\xi d\eta d\zeta \text{ for any }\hat{p}\in \mathcal{J}_{k-1}(\eta)$,\\

\item
$\int_{\hat{K}} \hat{q}_3 \hat{p}d\xi d\eta d\zeta \text{ for any }\hat{p}\in \mathcal{J}_{k-1}(\zeta)$,\\
\item
$\int_{\hat{K}} \hat{q}\cdot \hat{p}d\xi d\eta d \zeta \text{ for any }\hat{p}\in \Psi_{k-1}(\hat{K})$.
\end{enumerate}
\end{lemma}
\begin{proof}  Since the  dimensions of  the space $H_k(\hat{K})$ is equal to the number of these conditions,  it suffices to
 prove that $\hat{q}\equiv 0$ if  these conditions vanish. Since $\hat{q}\cdot \hat{\nu}\in Q_{k-1}(\hat{e})$, the first  condition (1)
  implies that
 \begin{equation*}
 \begin{split}
  \hat{q}_1=(1-\xi^2)(\hat{g}_1+\hat{f}_1),\\[0.5ex]
  \hat{q_2}=(1-\eta^2)(\hat{g}_2+\hat{f}_2),\\[0.5ex]
    \hat{q_2}=(1-\zeta^2)(\hat{g}_3+\hat{f}_3),\\[0.5ex]
  \end{split}
  \end{equation*}
  where $(\hat{g}_1\,, \hat{g}_2, \hat{g}_3)^T\in \Psi_{k-1}(\hat{K})$,
$\hat{f}_1\in \mathcal{J}_{k-1}(\xi)$,
$\hat{f}_2\in \mathcal{J}_{k-1}(\eta)$,
and $\hat{f}_3\in \mathcal{J}_{k-1}(\zeta)$.  Note that
  \begin{equation*}
  \int_{\hat{K}}(1-\xi^2)\hat{g}_1\hat{p}d\xi d\eta d\zeta=0 \text{ for any }\hat{p}\in \mathcal{J}_{k-1}(\xi).
  \end{equation*}
  By the condition (2),  this shows $\hat{f}_1=0$.   A similar argument (using the conditions (3)--(4))  yields
  $$
  \hat{f}_2=\hat{f}_3=0.
  $$
  Finally the condition (5) proves $\hat{g}_1=\hat{g}_2=\hat{g}_3=0$, which completes the proof.
  \end{proof}

To  design  finite element spaces for the components $\sigma_{12}$, $\sigma_{13}$, and $\sigma_{23}$,  of the shear stress, we need the following space
$$
S_k(X,Y)\times P_{k-1}(Z)\text{ for any }(X,Y,Z)\in \hat{K}:=[-1,1]^3.
$$
\begin{lemma}\label{lemma5.2}
Given any $\tau_{12}\in S_k(X,Y)\times P_{k-1}(Z)$, it can be uniquely determined
by the following conditions:
  \begin{enumerate}
  \item  the values of $\tau_{12}$ at $k$ distinct points on each edge of $\hat{K}$ that is  perpendicular to the  $(X,Y)$-plane,\\
\item the values of $\tau_{12}$ at $k(k-1)$ distinct points in the interior of each face of
$\hat{K}$ that parallels to the $Z$-axis,\\
  \item  the moments $\int_{\hat{K}} \tau_{12} p_{k-4}dXdYdZ$ for any $p_{k-4}\in P_{k-4}(X,Y)\times P_{k-1}(Z)$.
  \end{enumerate}
  Here the points in the second term are chosen in this way so that they lie in the sam plane as the points in the first term.
\end{lemma}
\begin{proof} Since $S_k(X,Y)$ is the space of the serendipity element of order $k$
with  respect to the variables $X$ and $Y$,  on each rectangle that
 parallels to the $(X,Y)$-plane,  $\tau_{12}$ can be uniquely determined by
\begin{equation*}
\begin{split}
 & \text{ the values of $\tau_{12}$ at four vertices of the rectangle},\\
&\text{ the values of $\tau_{12}$ at $k-1$ distinct points in the interior of each edge of
 the rectangle},\\
&\text{  the moments of  order $k-4$ of $\tau_{12}$ on the rectangle}.
\end{split}
  \end{equation*}
  Then the desired result follows from the fact that $S_k(X,Y)\times P_{k-1}(Z)$ is a product space.
\end{proof}

Then, on element $K$, the space for the stress can be defined as
\begin{equation}
\begin{split}
\Sigma_k(K):=\{\sigma\in H(\div, K, \S)|\sigma_n\in H_k(K), \sigma_{12}\in S_k(x, y)\times P_{k-1}(z),\\[0.5ex]
 \sigma_{13}\in S_k(x,z)\times P_{k-1}(y), \sigma_{23}\in S_k(y,z)\times P_{k-1}(x)\}.
\end{split}
\end{equation}

The global space is defined as
\[
\Sigma_k(\cT_h):=\{\tau\in \Sigma, \tau|_K\in \Sigma_k(K) \text{ for any }K\in\cT_h\}.
\]
On each element $K$, the space for the displacement is taken as
$$
V_k(K):=(Q_{k-1}(K))^3\oplus\{(Q_{k,x}, 0, 0)^T\}\oplus\{(0, Q_{k, y},0)^T\}\oplus\{(0, 0,  Q_{k, z})^T\},
$$
where
\begin{equation*}
\begin{split}
Q_{k,x}=x^k(P_{k-1}(y)+P_{k-1}(z)),\\[0.5ex]
Q_{k,y}=y^k(P_{k-1}(z)+P_{k-1}(x)),\\[0.5ex]
Q_{k,z}=z^k(P_{k-1}(x)+P_{k-1}(y)).
\end{split}
\end{equation*}

Then the global space for the displacement  reads
\[
V_{k}(\cT_h):=\{v\in V, v|_K\in V_k(K) \text{ for any } K\in \cT_h\}.
\]
\begin{remark} The lowest order element (k=1) of this family has
 21 stress  and  6 displacement degrees of freedom  per element, which is the three dimensional element of \cite{HuManZhang2013}, see some degrees of freedom in Figure \ref{3D}.
\end{remark}

\def\boa{\begin{picture}(120,125)(  0.,  0.)
   \def\la{\vrule width.4pt height.4pt}
 \multiput(   0.67,   0.50)(   2.667,   2.000){ 13}
       {\multiput(0,0)(  0.1333,  0.1000){ 10}{\la}}
 \multiput(  36.00, 106.17)(   0.000,  -3.333){ 24}
       {\multiput(0,0)(  0.0000, -0.1667){ 10}{\la}}
 \multiput( 115.17,  27.00)(  -3.333,   0.000){ 24}
       {\multiput(0,0)( -0.1667,  0.0000){ 10}{\la}}
  \put(0,0){\line(1,0){80}}\put(0,80){\line(1,0){80}}
  \put(80,0){\line(4,3){36}} \put(80,80){\line(4,3){36}}
  \put(0,80){\line(4,3){36}}
 \put(0,0){\line(0,1){80}}\put(80,0){\line(0,1){80}}
  \put(116,107){\line(0,-1){80}}\put(116,107){\line(-1,0){80}}
  \end{picture}}
 \def\bo{\begin{picture}(120,120)(  0.,  0.)
   \def\la{\vrule width.4pt height.4pt}
 \multiput(   0.00,   0.00)(   2.667,   2.000){ 13}{\la}
 \multiput(  36.00, 107.00)(   0.000,  -3.333){ 24}{\la}
 \multiput( 116.00,  27.00)(  -3.333,   0.000){ 24}{\la}
  \put(0,0){\line(1,0){80}}\put(0,80){\line(1,0){80}}
  \put(80,0){\line(4,3){36}} \put(80,80){\line(4,3){36}}
  \put(0,80){\line(4,3){36}}
 \put(0,0){\line(0,1){80}}\put(80,0){\line(0,1){80}}
  \put(116,107){\line(0,-1){80}}\put(116,107){\line(-1,0){80}}
  \end{picture}}

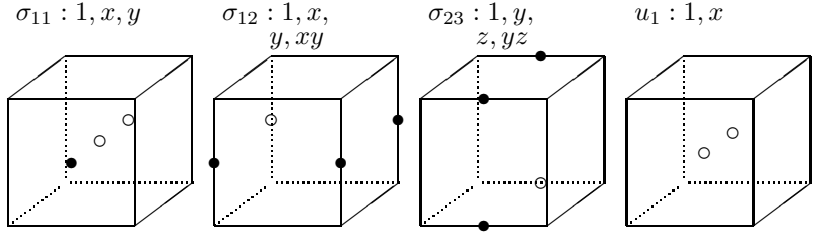
\begin{figure}[htb]
   \begin{center}\setlength{\unitlength}{0.6pt}
    \begin{picture}(520,135)(0,5)
 \put(0,0){\begin{picture}(100,100)(0,0) \put(0,0){\bo}
     \put(40,40){\circle*{6}} \put(58,53.5){\circle{6}}
         \put(76,67){\circle{6}}
   \put(5,130){$\sigma_{11}: 1, x, y$} \end{picture}}
  \put(130,0){\begin{picture}(100,100)(0,0) \put(0,0){\bo}
    \put(0,40){\circle*{6}}\put(80,40){\circle*{6}}
  \put(5,130){$\sigma_{12}: 1, x, $}
  \put(35,115){$ y,xy$}
       \put(116,67){\circle*{6}} \put(36,67){\circle{6}} \end{picture}}
  \put(260,0){\begin{picture}(100,100)(0,0) \put(0,0){\bo}
    \put(40,0){\circle*{6}}\put(40,80){\circle*{6}}
   \put(5,130){$\sigma_{23}: 1, y, $}
   \put(35,115){$ z, yz$}
       \put(76,107){\circle*{6}} \put(76,27){\circle{6}} \end{picture}}
    \put(390,0){\begin{picture}(100,100)(0,0) \put(0,0){\bo}
    \put(49,46.25){\circle{6}}\put(67,58.75){\circle{6}}
    \put(5,130){$u_1:1,x$}
         \end{picture}}
  \end{picture}
   \end{center}
\caption{Some nodal degrees of freedom. }\label{3D}
\end{figure}

To discretize the pure traction boundary  problem, we introduce the rigid motion space
$$
\text{RM}:=\text{span}\bigg\{
\begin{pmatrix} 1\\ 0\\ 0 \\ \end{pmatrix}, \begin{pmatrix} 0\\ 1 \\ 0\end{pmatrix}, \begin{pmatrix}0\\ 0 \\ 1\end{pmatrix},
 \begin{pmatrix}-y\\ x \\ 0\end{pmatrix}, \begin{pmatrix}-z\\ 0 \\ x\end{pmatrix}, \begin{pmatrix}0\\ -z \\ y\end{pmatrix}\bigg\},
$$
which defines
 \begin{equation}
   \begin{split}
    \Sigma_{k, 0}(\cT_h)
     &=\{ \tau\in  \Sigma_k(\cT_h) \mid
      \tau \nu=0 \ \quad  \text{ on  } \partial \Om\}, \\
	 V_{k,0}(\cT_h)
     &=\{ v\in V_k(\cT_h) \mid
       \ (v, w)_{L^2(\Omega)}=0 \quad \text{ for all } w\in \text{RM}
       \}.
       \end{split}
     	\end{equation}

It follows from the definitions of  the spaces $\Sigma_k(\cT_h)$ (resp. $\Sigma_{k, 0}(\cT_h)$) and
$V_k(\cT_h)$ (resp. $V_{k, 0}(\cT_h)$) that $\div \Sigma_k(\cT_h)\subset V_k(\cT_h)$ (resp. $\div \Sigma_{k, 0}(\cT_h)\subset V_{k, 0}(\cT_h)$).  Similar arguments of Theorems \ref{Theorem3.1}, \ref{Theorem4.1} and  \ref{Theorem4.2} can prove  the converses $V_k(\cT_h)\subset
\div \Sigma_k(\cT_h)$ and  $V_{k, 0}(\cT_h)\subset
\div \Sigma_{k, 0}(\cT_h)$, respectively.  In fact, to extend the result of Theorem \ref{Theorem3.1} to the present case, we only need  essentially
 three one-dimension-arguments used in Theorem \ref{Theorem3.1};  while to get a generalization of Theorems \ref{Theorem4.1} and  \ref{Theorem4.2}, we only need essentially three two-dimension-arguments used in Theorems \ref{Theorem4.1} and  \ref{Theorem4.2}.  In particular,  in this way, we can get three two-dimension-rigid motion spaces, which proves, on a  macroelement consisting of eight elements,  the kernel space $N_M$ (see \eqref{kernel} for the definition in two dimensions)   is  the rigid motion space  in three dimension. This in turn implies a similar result of Lemma \ref{Lemma4.3}. Finally,  this indicates the well-posedness of  this
family of elements.

\section{Reduced elements in both two and three dimensions}
In this section we present a  family of reduced elements for these in Sections 2 and 5. To this end, we introduce
Airy's stress function for a scalar field $q(X,Y)$ as follows
\begin{equation*}
J_{X,Y}(q(X,Y)):=
\begin{pmatrix}
\frac{\pa^2 q}{\pa Y^2}& -\frac{\pa^2 q}{\pa X\pa Y}\\
-\frac{\pa^2 q}{\pa X\pa Y} & \frac{\pa^2 q}{\pa X^2}
\end{pmatrix}.
\end{equation*}
Throughout this section we let $(X,Y,Z)$ denote  permutations of $(x,y,z)$.
\subsection{The reduced elements in two dimensions}
We define the shape function space for the BDFM element \cite{Brezzi-Douglas-Fortin-Marini} as
\[
BDFM_k(K):=(P_k(K))^2\backslash\ \sspan\{(0,x^k),(y^k, 0)\}.
\]
The stress space  of the reduced element of order $k$ is defined as
\[
\Sigma_k^{R}(K):=\{\tau\in\S, \tau_n\in BDFM_k(K), \tau_{12}\in P_k(x,y)\}\oplus E_k(K),
\]
where
$$
E_k(K):=\sspan\{J_{x,y}(x^{k+1}y^2), J_{x,y}(x^2y^{k+1})\}.
$$
The degrees of freedom for the stress are inherited from the  BDFM element and the
serendipity element:
\begin{enumerate}
\item  the moments of degree not greater than $k-1$ on the four edges of $K$
       for $\sigma_n\cdot \nu $,\\
\item the moments of  degree not greater than $k-2$ on  $K$ for $\sigma_n$,\\

\item  the values of $\sigma_{12}$ at four vertices of $K$,\\
\item the values of $\sigma_{12}$ at $k-1$ distinct points in the interior of each edge of
$K$,\\
\item the moments of degree not greater than $k-4$ on  $K$ for
$\sigma_{12}$.
\end{enumerate}
The global space for the stress of order $k$ is defined as
\begin{equation}
\Sigma_k^{R}(\cT_h):=\{\tau\in \Sigma, \tau|_K\in \Sigma_k^{R}(K) \text{ for any }K\in\cT_h\}.
\end{equation}
On each element $K$, the space for the displacement is taken as
$$
V_k^R(K):=(P_{k-1}(K))^2.
$$
Then the global space for the displacement  reads
\begin{equation}
V_{k}^R(\cT_h):=\{v\in V, v|_K\in V_k^R(K) \text{ for any } K\in \cT_h\}.
\end{equation}
\begin{remark} Let $RT_k(K)=RT_{k}(K)=P_{k, k-1}(K)\times P_{k-1, k}(K)$.  One can also define the space for the stress
 as
\[
\hat{\Sigma}_k(K):=\{\tau\in\S, \tau_n\in RT_k(K), \tau_{12}\in P_k(x,y)\}\oplus E_k(K).
\]
The space for the displacement in this case is
$$
\hat{V}_k(K):=(Q_{k-1}(K))^2.
$$
\end{remark}

\subsection{The reduced elements in three dimensions}
 On element $K\in\cT_h$,  for $k\geq 1$,  we define the  Raviart--Thomas element space by
 $$
 RT_{k}(K)=P_{k, k-1, k-1}(K)\times P_{k-1, k, k-1}(K)\times P_{k-1, k-1, k}(K),
 $$
 where
\begin{equation*}
\begin{split}
P_{k, k-1, k-1}(K)=P_k(x)\times P_{k-1}(y)\times P_{k-1}(z),\\[0.5ex]
P_{k-1, k, k-1}(K)=P_{k-1}(x)\times P_k(y)\times P_{k-1}(z),\\[0.5ex]
P_{k-1, k-1, k}(K)=P_{k-1}(x)\times P_{k-1}(y)\times P_k(z).
\end{split}
\end{equation*}
Given a scalar field $q(X,Y)$ and the corresponding Airy's function $J_{X,Y}(q(X,Y))$, we define $\tau(J_{X,Y}q(X,Y))\in H(\div, K, \S)$ such that
$$
\tau_{X,X}=(J_{X,Y}q(X,Y))_{X,X}, \tau_{Y,X}=\tau_{X,Y}=(J_{X,Y}q(X,Y))_{X, Y}, \tau_{Y,Y}=(J_{X,Y}q(X,Y))_{Y, Y}
$$
and the rest entries are zero. This notation allows to define
\begin{equation*}
\begin{split}
E_{k}(K):&=\sspan\bigg\{\tau(J_{x,y}(x^{k+1}y^2)), \tau(J_{x,y}(x^2y^{k+1})) \bigg\}P_{k-1}(z)\\
&\oplus \sspan\bigg\{\tau(J_{x,z}(x^{k+1}z^2)), \tau(J_{x,z}(x^2z^{k+1}))\bigg\}P_{k-1}(y)\\
&\oplus \sspan\bigg\{\tau(J_{y,z}(y^{k+1}z^2)), \tau(J_{y,z}(y^2z^{k+1}))\bigg\}P_{k-1}(x).
\end{split}
\end{equation*}
Then, on element $K$, the space for the stress can be defined as
\begin{equation}
\begin{split}
\Sigma_k^R(K):=\{\sigma\in H(\div, K, \S)|\sigma_n\in RT_k(K), \sigma_{12}\in P_k(x, y)\times P_{k-1}(z),\\[0.5ex]
 \sigma_{13}\in P_k(x,z)\times P_{k-1}(y), \sigma_{23}\in P_k(y,z)\times P_{k-1}(x)\}\oplus E_k(K).
\end{split}
\end{equation}
 The stress $\tau\in \Sigma_k^R(K)$  can be uniquely
determined by the following conditions:
\begin{enumerate}
\item
$\int_e \tau_n\cdot \nu p ds \text{ for any }p \in Q_{k-1}(e)
\text{ and  any }e \subset \partial K$,\\
\item
$\int_K \tau_n \cdot pdx dy dz  \text{ for any }p\in \Psi_{k-1}(K)$,\\
  \item  the values of $\tau_{XY}$ at $k$ distinct points on each edge of $K$ that is  perpendicular to the  $(X,Y)$-plane,\\
\item the values of $\tau_{XY}$ at $k(k-1)$ distinct points in the interior of each face of
$K$ that parallels to the $Z$-axis,\\
  \item  the moments $\int_K \tau_{XY} p_{k-4}dXdYdZ$ for any $p_{k-4}\in P_{k-4}(X,Y)\times P_{k-1}(Z)$.
\end{enumerate}
The proof for unisolvence of these degrees of freedom  follows directly from those  the  RT element and the
serendipity element, which is omitten herein; c.f. similar proofs in Lemmas \ref{lemma5.1} and \ref{lemma5.2}.

The global space is defined as
\[
\Sigma_k^R(\cT_h):=\{\tau\in \Sigma, \tau|_K\in \Sigma_k^R(K) \text{ for any }K\in\cT_h\}.
\]
On each element $K$, the space for the displacement is taken as
$$
V_k^R(K):=(Q_{k-1}(K))^3.
$$
Then the global space for the displacement  reads
\[
V_{k}^R(\cT_h):=\{v\in V, v|_K\in V_k^R(K) \text{ for any } K\in \cT_h\}.
\]

To discretize the pure traction boundary  problem, we  define
 \begin{equation}
   \begin{split}
    \Sigma_{k, 0}^R(\cT_h)
     &=\{ \tau\in  \Sigma_k^R(\cT_h) \mid
      \tau \nu=0 \ \quad  \text{ on  } \partial \Om\}, \\
	 V_{k,0}^R(\cT_h)
     &=\{ v\in V_k^R(\cT_h) \mid
       \ (v, w)=0 \quad \text{ for all } w\in \text{RM}
       \}.
       \end{split}
     	\end{equation}

It follows from the definitions of  the spaces $\Sigma_k^R(\cT_h)$ (resp. $\Sigma_{k, 0}^R(\cT_h)$) and
$V_k^R(\cT_h)$ (resp. $V_{k, 0}^R(\cT_h)$) that $\div \Sigma_k^R(\cT_h)\subset V_k^R(\cT_h)$ (resp. $\div \Sigma_{k, 0}^R(\cT_h)\subset V_{k, 0}^R(\cT_h)$).  Similar arguments of Theorems \ref{Theorem3.1}, \ref{Theorem4.1} and  \ref{Theorem4.2} can prove  the converses $V_k^R(\cT_h)\subset
\div \Sigma_k^R(\cT_h)$ and  $V_{k, 0}^R(\cT_h)\subset
\div \Sigma_{k, 0}^R(\cT_h)$, respectively. This indicates the well-posedness of  this
family of elements.

\begin{remark} The lowest order element (k=1) of this family has
 8 stress and 2 displacement, and  18 stress  and  3 displacement degrees of freedom  per element for two and three dimensions, respectively,
  which were announced independently  by Chen and his collaborators \cite{Chen2013} after the first version of the paper was submitted.
\end{remark}

\section{The error estimate and numerical results}
\subsection{The error estimate}
The section is devoted to the error analysis of the approximation
defined by \eqref{discrete}.  It follows from \eqref{continuous} and
\eqref{discrete} that
\begin{equation}\label{error}
\begin{split}
&(A(\sigma-\sigma_{k,h}),\tau_h)_{L^2(\Om)}+(\div\tau_h,(u-u_{k,h}))_{L^2(\Om)}=0
\text{ for any }\tau_h\in\Sigma_k(\cT_h),\\
&(\div(\sigma-\sigma_{k,h}),v_h)_{L^2(\Om)}=0 \text{ for any }v_h\in V_k(\cT_h).\\
\end{split}
\end{equation}
Let $P_h$ be the $L^2$ projection operator from $L^2(\Omega, \R^n)$
onto $V_k(\cT_h)$.  Since $\div\sigma_{k,h}\in V_k(\cT_h)$, the
second equation of \eqref{error} yields
\begin{equation}\label{eq4.2}
\begin{split}
\|\div(\sigma-\sigma_{k,h})\|_{L^2(\Omega)}&=\|\div\sigma-P_h\div\sigma\|_{L^2(\Omega)}\\
&\leq Ch^{m} |\div\sigma |_{H^m(\Omega)}\text{ for any } 0\leq m\leq k.
\end{split}
\end{equation}
It follows from the K-ellipticity,  Theorem \ref{Theorem3.1}
and the approximation properties of $\Sigma_{k,h}$ and $V_{k,h}$ that
\begin{equation}
\begin{split}
&\|\sigma-\sigma_{k,h}\|_{L^2(\Omega)}+\|u-u_{k,h}\|_{L^2(\Omega)}\\
&\leq
C\big(\inf\limits_{\tau_h\in\Sigma_{k,
h}}\|\sigma-\tau_h\|_{H(\div,\Omega)}+\inf\limits_{v_h\in V_{k,h}}\|u-v_h\|_{L^2(\Omega)}\big )\\
&\leq Ch^{m} (|\sigma|_{H^{m+1}(\Omega)}+|u|_{H^m(\Omega)}) \text{ for any } 0\leq m\leq k.
\end{split}
\end{equation}
A similar error estimate holds for the pure traction boundary  problem studied in section 4 and the reduced elements in Section 6.
\begin{remark} By using the mesh  dependent norm in Subsection 4.2, cf. \eqref{eq4.23} and \eqref{eq4.25},  we can  get an improved error estimate:
$$
\|\sigma-\sigma_{k,h}\|_{L^2(\Omega)}\leq Ch^m|\sigma|_{H^{m}(\Omega)} \text{ for any } 0\leq m\leq k.
$$
\end{remark}
\subsection{The numerical result}
The first example is  presented to demonstrate the second order method (with $k=2$)
for the pure displacement boundary  problem with a homogeneous boundary condition
         that $u\equiv 0$ on $\partial\Omega$; see \cite{HuManZhang2013} for numerical examples for $k=1$.
          Assume the material is isotropic in the sense that
 \begin{equation*}
A \sigma = \frac{1}{2\mu} \left(
	   \sigma - \frac{\lambda}{2\mu + 2 \lambda} \operatorname{tr}(\sigma)
		\delta \right),
\end{equation*}
where $\delta$ is the identity matrix, and $\mu$ and $\lambda$ are the
Lam\'e constants such that $0<\mu_1\leq \mu \leq \mu_2$ and
$0<\lambda<\infty$.  In the numerical example, these parameters are chosen as
$$
\lambda=1\quad  \mu =\frac{1}{2}.
$$

Let the  exact solution on the unit square $[0,1]^2$ be
 \begin{equation}
 u=(\sin\pi x\sin \pi y, \sin \pi x \sin \pi y)^T.
 \end{equation}

 \begin{table}[htb]
  \caption{ The error and the order of convergence.}\label{Table1}
\begin{center}  \begin{tabular}{c|cr|cr|cr}  
\hline & $ \| u- u_{2,h}\|_{0}$ & rate &
    $ \|\sigma -\sigma_{2,h}\|_{0}$ & rate  &
    $ \|\div(\sigma -\sigma_{2,h})\|_{0}$ & rate  \\ \hline
 1& 0.3156  &0.0&  2.0116&0.0&   7.8083&0.0\\
 2&   0.0693&2.2&  0.4465&2.2&   1.9752&2.0\\
 3&  0.0166&2.1&    0.1134&2.0& 0.4760&2.1\\
 4&  0.0041&2.0&  0.0285&2.0&   0.1175&2.0\\
 5&0.0010&2.0& 0.0071&2.0&    0.0293&2.0\\
 6& 2.5408e-004&2.0&  0.0018&1.9&  0.0073&2.0\\
 7&  6.3503e-005&2.0&  4.4605e-004&2.0&   0.0018&2.0 \\
      \hline
\end{tabular}\end{center} \end{table}

In the computation, the level one grid is the given domain, a unit
   square or a unit cube.
Each grid is refined into a half-size grid uniformly, to get
   a higher level grid, see the first column in Table \ref{Table1}.

As the second example,  we compute the pure traction boundary  problem
 with the exact solution
\begin{equation}\label{s-l} u= \left[100 x^2 (1-x)^2 y^2 (1-y)^2 -\frac 19 \right]
    \begin{pmatrix} 1\\-1\end{pmatrix}.
    \end{equation}
The matrix $A$ is same as that in the first  example.
Our new finite element has no problem in solving the pure
  traction boundary  problems.
The convergence results are listed in Table \ref{t-s-l}.

 \begin{table}[htb]
  \caption{\label{t-s-l}
 The errors and the order of convergence for
   the pure traction boundary  problem}
\begin{center}  \begin{tabular}{c|cc|cc|cc}  
\hline & $ \| u- u_{2, h}\|_{0}$ & rate&
    $ \|\sigma -\sigma_{2, h}\|_{0}$ & rate  &
    $ \|\div(\sigma -\sigma_{2, h})\|_{0}$ & rate  \\ \hline
 2&  0.0264&0.0&  0.2516&0.0&   2.4645&0.0 \\
 3&  0.0107&1.3&  0.0804&1.6&   0.7090&1.8 \\
 4&  0.0029&1.9&  0.0211&2.0&   0.1807&2.0 \\
 5&  7.2940e-004&2.0&   0.0054&2.0&   0.0453&2.0 \\
 6&  1.8315e-004&2.0& 0.0013&2.0&   0.0113&2.0 \\
 7&  4.5836e-005&2.0&  3.3684e-004&2.0&   0.0028&2.0 \\
      \hline
\end{tabular}\end{center} \end{table}


\begin{thebibliography}{999}
\bibitem{Adams-Cockburn} S. Adams and B. Cockburn,
  A mixed finite element method for elasticity in three dimensions, J.
    Sci. Comput. 25 (2005), no. 3, 515--521.

\bibitem{Amara-Thomas} M. Amara and J. M. Thomas,
   Equilibrium finite elements for the linear elastic problem, Numer.
    Math. 33 (1979), 367--383.

\bibitem{Arnold-Awanou} D. N. Arnold and G. Awanou,
   Rectangular mixed finite elements for elasticity, Math. Models
    Methods Appl. Sci. 15 (2005),  1417--1429.

\bibitem{Arnold-Awanou2011} D. N. Arnold and G. Awanou, The serendipity family of finite elements,
Found. Comput. Math.
 11(2011), 337--344.


\bibitem{Arnold-Awanou-Winther}
    D. Arnold, G. Awanou and R. Winther,
   Finite elements for symmetric tensors in three dimensions,
    Math. Comp. 77 (2008), no. 263, 1229--1251.


\bibitem{Arnold-Brezzi-Douglas}
     D. N. Arnold, F. Brezzi and J. Douglas, Jr.,
     PEERS: A new mixed finite element for plane elasticity,
    Jpn. J. Appl. Math. 1 (1984),  347--367.

\bibitem{Arnold-Douglas-Gupta}
  D. N. Arnold, J. Douglas Jr., and C. P. Gupta,
    A family of higher order mixed finite element
    methods for plane elasticity, Numer. Math. 45 (1984),  1--22.


\bibitem{Arnold-Falk-Winther}
  D.N. Arnold, R. Falk and R. Winther,
  Mixed finite element methods for linear elasticity with
    weakly imposed symmetry, Math. Comp. 76 (2007), no. 260, 1699--1723.


\bibitem{Arnold-Winther-conforming}
    D. N. Arnold and R. Winther,
    Mixed finite element for elasticity, Numer. Math. 92 (2002),
     401--419.

\bibitem{Arnold-Winther-n} D. N. Arnold and R. Winther,
   Nonconforming mixed elements for elasticity, Math. Models.
   Methods Appl. Sci. 13 (2003), 295--307.

\bibitem{Awanou} G. Awanou,
    Two remarks on rectangular mixed finite elements
   for elasticity, J. Sci. Comput. 50 (2012), 91--102.

 \bibitem{BecacheJolyTsogka2002}E. B\'{e}cache, P. Joly and C. Tsogka, A new family of mixed finite elements for
 the linear elastodynamic problem ,  SIAM J. Numer. Anal., 39(2002), pp. 2109--2132.


\bibitem{Boffi-Brezzi-Fortin}
   D. Boffi, F. Brezzi and M. Fortin,
   Reduced symmetry elements in linear elasticity, Commun.
    Pure Appl. Anal. 8 (2009), no. 1, 95--121.

\bibitem{BrennerScott} S.~C.~Brenner and L.~ R.~ Scott, The
mathematical theorey of finite element methods, Springer-Verlag,
1996.


\bibitem{Brezzi-Douglas-Fortin-Marini}F. Brezzi, J. Douglas, Jr., M.
Fortin,  L. D. Marini,  Efficient rectangular mixed finite elements
in two and three space variables, A. I. R.O., Mode. Math. Anal.
Numer., 21(1987): 581--604.


\bibitem{Brezzi-Fortin} F. Brezzi and  M. Fortin,
   Mixed and hybrid finite element methods, Springer, 1991.

\bibitem{Carstensen} C. Carstensen, M. Eigel and J. Gedicke,
   Computational competition of symmetric mixed FEM in linear elasticity,
   Comput. Methods Appl. Mech. Engrg. 200 (2011), no. 41-44, 2903--2915.

\bibitem{Carstensen-Gunther-Reininghaus-Thiele2008}C. Carstensen,  D. G\"{u}nther, J. Reininghaus, J. Thiele,
The Arnold--Winther mixed FEM in linear elasticity. Part I: Implementation and numerical verification,
Comput. Methods Appl. Mech. Engrg. 197 (2008) 3014--3023.


\bibitem{Chen-Wang} S.-C. Chen and Y.-N. Wang,
   Conforming rectangular mixed finite elements for elasticity,
  J. Sci. Comput. 47 (2011), no. 1, 93--108.

\bibitem{Chen2013} S. C. Chen,  Presentation in the workshop on ``Finite element methods and its applications",  Beijing, China,  December 7 2013.

\bibitem{CiaBook} P. G.  Ciarlet,  The finite element method
for elliptic problems,  North--Holland, 1978;
        reprinted as SIAM Classics in Applied Mathematics, 2002.


\bibitem{Cockburn-Gopalakrishnan-Guzman}
  B. Cockburn, J. Gopalakrishnan and J. Guzm\'an,
    A new elasticity element made for enforcing weak stress symmetry,
    Math. Comp. 79 (2010), no. 271, 1331--1349.



\bibitem{GiraultRaviart1986} V. Girault and P. A. Raviart, Finite
element methods for Navier-Stokes equations: theory and algorithms,
Springer-Verlag, Berlin, Heidelberg  1986.

\bibitem{Gopalakrishnan-Guzman-n} J. Gopalakrishnan and J. Guzm\'an,
  Symmetric nonconforming mixed finite elements for linear elasticity,
 SIAM J. Numer. Anal. 49 (2011), no. 4, 1504--1520.


\bibitem{Gopalakrishnan-Guzman} J. Gopalakrishnan and J. Guzm\'an,
    A second elasticity element using the matrix bubble,
    IMA J. Numer. Anal. 32 (2012), no. 1, 352--372.


\bibitem{Guzman} J. Guzm\'an,
   A unified analysis of several mixed methods for elasticity
   with weak stress symmetry,
    J. Sci. Comput. 44 (2010), no. 2, 156--169.

\bibitem{HW} Jason S. Howell,  Noel J. Walkington,
 Inf-sup conditions for twofold saddle point problems.
 Numer. Math. 118(2011), No. 4, 663--693.

\bibitem{HuManZhang2013}J. Hu, H. Y. Man, and S. Y. Zhang, A simple conforming mixed finite
   element for linear elasticity on rectangular grids  in any space dimension, J. Sci. Comput.
DOI 10.1007/s10915-013-9736-6.

\bibitem{HuManZhang2012}J. Hu, H. Y. Man, and S. Y. Zhang,
A minimal mixed finite element method
   for linear elasticity in the symmetric formulation
   on $n$-rectangular grids, arXiv:1304.5428[math.NA] (2013).


\bibitem{Hu-Shi}  J. Hu and Z. C. Shi,
     Lower order rectangular nonconforming mixed elements for plane elasticity,
     SIAM J. Numer. Anal.  46 (2007),   88--102.

\bibitem{Johnson-Mercier} C. Johnson and B. Mercier,
    Some equilibrium finite element methods for two-dimensional elasticity
     problems, Numer.Math. 30 (1978),  103--116.

\bibitem{Man-Hu-Shi} H.-Y. Man, J. Hu and Z.-C. Shi,
   Lower order rectangular nonconforming mixed finite element for
     the three-dimensional elasticity problem,
   Math. Models Methods Appl. Sci. 19 (2009), no. 1, 51--65.

\bibitem{Morley} M. Morley, A family of mixed finite elements
   for linear elasticity Numer. Math. 55 (1989), no. 6, 633­-666.


\bibitem{Stenberg-1}
     R. Stenberg, On the construction of optimal mixed finite
     element methods for the linear elasticity
     problem, Numer. Math. 48 (1986),   447--462.



\bibitem{Stenberg-2}
     R. Stenberg, Two low-order mixed methods for the elasticity problem,
     In: J. R. Whiteman (ed.):
     The Mathematics of Finite Elements and Applications, VI. London:
     Academic Press, 1988,   271--280.

\bibitem{Stenberg-3} R. Stenberg,
   A family of mixed finite elements for the elasticity problem,
    Numer. Math. 53 (1988),  no. 5, 513--538.

\bibitem{Stenberg90} R. Stenberg,  A technique for analysing finite element methods for viscous incompressible flow,
Internat. J. Numer. Methods Fluids, 11 (1990), pp. 935--948.

\bibitem{Yi} S. Y. Yi,
    A New nonconforming mixed finite element method for linear elasticity,
    Math. Models  Methods Appl. Sci. 16 (2006),  979--999.

\end{thebibliography}
\end{document}